%% file: pevp2.tex
\documentclass[final]{siamltex}
\usepackage{showkeys}
\usepackage{amssymb,amsmath}
\usepackage{
  tikz,%
  pgfplots, %
  xcolor, %
  mathdots, %
}
\pgfplotsset{compat=1.15}
\usetikzlibrary{calc}
\usetikzlibrary{arrows}
\usetikzlibrary{patterns}
\usetikzlibrary{decorations.pathmorphing}
\usepgfplotslibrary{statistics}

\pgfplotscreateplotcyclelist{warmercolors}{
  {SPECblue!80!SPECorange},
  {SPECblue!70!SPECorange},
  {SPECblue!60!SPECorange},
  {SPECblue!50!SPECorange},
  {SPECblue!40!SPECorange},
  {SPECblue!30!SPECorange},
  {SPECblue!20!SPECorange},
  {SPECblue!10!SPECorange},
  {SPECorange},
  {SPECorange!90!SPECred},
  {SPECorange!80!SPECred},
  {SPECorange!70!SPECred},
  {SPECorange!60!SPECred},
  {SPECorange!50!SPECred},
  {SPECorange!40!SPECred},
  {SPECorange!30!SPECred},
  {SPECorange!20!SPECred},
  {SPECorange!10!SPECred},
  {SPECred},
  {SPECred!90!black},
  {SPECred!80!black},
  {SPECred!70!black},
  {SPECred!60!black},
  {SPECred!50!black},
  {SPECred!40!black},
  {SPECred!30!black},
  {SPECred!20!black},
  {SPECred!10!black},
  {black},
}

\definecolor{SPECorange}{rgb}{1.0,.5625,0}
\definecolor{SPECblue}{rgb}{0,0,0.75}
\definecolor{SPECred}{rgb}{0.75,0,0}
\definecolor{SPECgreen}{rgb}{0,0.75,0}
\definecolor{SPECblack}{rgb}{0.75,0.75,0.75}
\definecolor{SPECCorange}{rgb}{1,0.9,0.6}
\definecolor{SPECCblue}{rgb}{0,0,0.25}
\definecolor{SPECCred}{rgb}{1.0,0,0}
\definecolor{SPECCgreen}{rgb}{0.25,0.75,0.25}

\usepackage{cite}
\usepackage{booktabs}
\usepackage{multirow}
\usepackage{listings}
\usepackage[ruled, lined, boxed, algo2e]{algorithm2e}

\newtheorem{remark}[theorem]{\it Remark}
\newtheorem{example}[theorem]{\it Example}

\usepackage{ifpdf}
\ifpdf
\usepackage[pdftex,colorlinks=true,linkcolor=SPECred,filecolor=SPECgreen,citecolor=SPECblue,pdfpagemode=UseNone,hypertexnames=false]{hyperref}
\else
\usepackage[dvips,colorlinks=true,linkcolor=SPECred,filecolor=SPECgreen,citecolor=SPECblue,pdfpagemode=UseNone]{hyperref}
\fi

\providecommand{\abs}[1]{\left\lvert#1\right\rvert} %
\newcommand{\norm}[1]{\mbox{$\parallel\!#1\!\parallel$}}

\newcommand{\R}{\ensuremath{\mathbb{R}}}

\newcommand{\Rnn}{\ensuremath{\mathbb{R}^{n\times n}}}

\newcommand{\cN}{\ensuremath{\mathcal{N}}}

\newcommand{\dmu}{\,\mathrm{d}\mu}

\title{Solving the Parametric Eigenvalue Problem by Taylor Series and Chebyshev Expansion%
  \thanks{The research has been partially funded by the Deutsche
Forschungsgemeinschaft (DFG)---Project-ID 318763901---SFB1294.}
}

\author{Thomas Mach\footnotemark[3]\and Melina A.\ Freitag\footnotemark[3]}

\begin{document}
\maketitle

\renewcommand{\thefootnote}{\fnsymbol{footnote}}

\footnotetext[3]{%
  University of Potsdam, Institute of Mathematics, 
  Karl-Liebknecht-Str.\ 24--25, 14476 Potsdam, Germany;
  \mbox{(\texttt{\{thomas.mach,melina.freitag\}@uni-potsdam.de})}.}%

\renewcommand{\thefootnote}{\arabic{footnote}}

\date{\today}
\maketitle

\begin{abstract}
  We discuss two approaches to solving the parametric (or stochastic) eigenvalue
  problem. One of them uses a Taylor expansion and the other a Chebyshev
  expansion. The parametric eigenvalue problem assumes that the matrix $A$
  depends on a parameter $\mu$, where $\mu$ might be a random
  variable. Consequently, the eigenvalues and eigenvectors are also functions of
  $\mu$. We compute a Taylor approximation of these functions about $\mu_{0}$ by
  iteratively computing the Taylor coefficients. The complexity of this approach
  is $O(n^{3})$ for all eigenpairs, if the derivatives of $A(\mu)$ at $\mu_{0}$
  are given.  The Chebyshev expansion works similarly. We first find an initial
  approximation iteratively which we then refine with Newton's method.  This
  second method is more expensive but provides a good approximation over the
  whole interval of the expansion instead around a single point.
  
  We present numerical experiments confirming the complexity and demonstrating
  that the approaches are capable of tracking eigenvalues at intersection
  points.  Further experiments shed light on the limitations of the Taylor
  expansion approach with respect to the distance from the expansion point
  $\mu_{0}$.
\end{abstract}



\pagestyle{myheadings} %
\thispagestyle{plain} %
\markboth{T.\ Mach and M. A. Freitag
}%
{Solving Parametric Eigenvalue Problems}

\section{Introduction}
\label{sec:introduction}

We are concerned with the \emph{parametric eigenvalue problem}
\begin{align}
  A(\mu) v(\mu) = v(\mu) \lambda(\mu), \label{eq:sevp}
\end{align} 
where the $n\times n$ entries of $A(\mu)$ are functions of $\mu$.
If $\mu$ is a random variable, this is referred to as stochastic eigenvalue
problem. A first example of such a problem is
\begin{align}
  A(\mu) = \exp_{\circ}(-\mu U), \label{eq:example-proof-of-concept}
\end{align}
where $\exp_{\circ}(\cdot)$ is the entrywise (or Hadamard) exponential---think
Matlab's \texttt{exp} instead of \texttt{expm}.\footnote{For \texttt{expm} the
  solution is known and can be found in many textbooks on numerical linear
  algebra. Since the eigenvectors do not depend on $\mu$ in the case of
  \texttt{expm} the problem simplifies to finding the eigenvalue decomposition
  of $U$ and applying \texttt{expm} to the Jordan form.} If the matrix $U$
contains in $U_{ij}$ the pointwise distances between points $x_{i}$ and $x_{j}$,
then this is a kernel method to analyze the structure of these point sets, see,
for instance, \cite{haasdonk2007invariant}. It is known that if $U$ is almost
negative definite and $\mu>0$, then $A$ is positive semidefinite
\cite[Thm.~6.3.6]{b172}. Thus, arguably, this is one of the simplest, yet
interesting parametric eigenvalue problems. The functions for the entries of
$A(\mu)$ are all of similar structure and are differentiable as often as
required. If, furthermore, all points $x_{i}$ are disjoint, the matrix is even positive
definite for all $\mu>0$. Thus, an algorithm failing for this example is likely
not useful. 

Other examples can be found if $A$ is the discretization of a partial
differential operator and $\mu$ an unknown parameter of the physical model, say
the mass of a passenger in a car, the temperature, or an unknown material
constant. Matrix-valued ODEs, for instance\begin{align*}
  \dot{C}(t) = -\alpha C(t) A^{T}\Gamma A C(t), \quad C(0) = C_{0}=C_{0}^{T}
\end{align*}
describing the dynamics of the covariance matrix in an ensemble Kalman
inversion, see \cite{bungert2022} for details, are another source of parametric
eigenvalue problems.

A few things about \eqref{eq:sevp} are obvious. If $(v(\mu),\lambda(\mu))$ is a
solution of \eqref{eq:sevp} then so is $(\gamma v(\mu),\lambda(\mu))$ for all
$\gamma\neq 0$, if normalization of the eigenvector is required for all $\gamma$
with $\abs{\gamma}=1$.  Fixing~$\mu$ turns \eqref{eq:sevp} into a standard
eigenvalue problem. For a non-defective $n\times n$ matrix, there are $n$
eigenpairs. If all the function $\left.A(\mu)\right|_{ij}$ are continuous in
$\mu$, then so are the $n$ functions $\lambda(\mu)$ and $v(\mu)$.\footnote{Note
  that there might be $\mu$ for which $A$ is defective. In this case the
  eigenvalues and eigenspaces are still continuous.} In fact, if $A(\mu)$
consists of holomorphic functions in $\mu$, then the eigenvalues $\lambda(\mu)$
are analytic functions with only algebraic singularities, see \cite[Ch. 2]{b473}
for details. Singularities do not occur if $A(\mu)$ is symmetric, see
\cite[Ch. 7]{b473}.  In this paper we will assume that the functions
$\left.A(\mu)\right|_{ij}$ are sufficiently smooth in $\mu$.

The problem \eqref{eq:sevp} has been investigated in the past. The following is
a selection---without claim of completeness---of the available literature on the
parametric and stochastic eigenvalue problem. Very recently Ruymbeek,
Meerbergen, and Michiels~\cite{ruymbeek2022tensor} used a tensorized Arnoldi
method to compute the extreme eigenvalues of \eqref{eq:sevp} using a polynomial
chaos expansion,
\begin{align*}
  A(\mu) = \sum_{\ell=0}^{M} A_{\ell}\phi_{\ell}(\mu),
\end{align*}
to discretize the problem in the parameter space. Soused\'ik and Elman
\cite{sousedik2016inverse} in a fairly similar approach used an inverse
iteration to find the smallest eigenvalue(s) of \eqref{eq:sevp} also using a
polynomial chaos expansion.  The resulting problem is a nonlinear tensor
equation, which is solved with a generalization of the inverse subspace
iteration. Ghanem and Ghosh \cite{ghanem2007efficient} used a similar idea,
but employed the Newton-Raphson algorithm to solve the tensor equation. Rahman
\cite{rahman2006solution} computed statistical moments of the generalized
stochastic eigenvalue problem,
\begin{align}
  A(\mu)v(\mu) = \lambda(\mu)B(\mu)v(\mu). \label{eq:sgevp}
\end{align}
Verhoosel, Guti\'errez, and Hulshoff \cite{verhoosel2006iterative} used a finite
element discretization for $\lambda(\mu)$ and $v(\mu)$ while assuming
$A(\mu) = A_{0} + \sum_{i=1}^{m}A_{i}\mu_{i}$, where $\mu$ is a parameter
vector, to solve a parametric eigenvalue problem. This discretization also leads
to a tensor structured equation, which they solved using an inverse iteration.  A
different approach is taken by Williams
\cite{williams2010method,williams2013method}, who adds an artificial time
dependency to \eqref{eq:sevp} to turn the problem into an integral
equation. Williams also uses a polynomial chaos expansion, however, for the
eigenvectors instead of $A(\mu)$.

All these methods have in common that they only compute a few of the smallest or
largest eigenvalues using inverse iteration, the power method, or more generally
a Krylov-type method. Some authors focus on the stochastic nature of the
problem. A different approach was taken very recently by Alghamdi, Boffi, and
Bonizzoni \cite{alghamdi2022greedy}. They are interested in a parameter
dependent PDE and use a model order reduction approach based on sparse grids to
find crossing points of the eigenvalue functions. They choose the parameter
$\mu$ from the $\R^{d}$.

In contrast we consider the parameter $\mu\in\R$ and focus on finding the
functions $\lambda(\mu)$ and $v(\mu)$ for the eigenvalues we are interested
in. If $A(\mu)$ is not too large, we can easily pick \emph{any}---not just the
smallest or largest---eigenvalue of $A(\mu_{0})$ and track this eigenvalue over
an interval or in the neighborhood of~$\mu_{0}$. The approach presented here is
not restricted to the smallest or largest eigenvalues. In fact we can compute
approximations for all eigenvalues, notwithstanding that in many application
only a few eigenvalues are desired.

Other generalizations of the standard eigenvalue problem,
\begin{align}
  Av = \lambda v,
  \label{eq:evp}
\end{align}
have been studied as well, far most the generalized eigenvalue problem
\cite{MolS73},
\begin{align}
  Av = \lambda Bv.
  \label{eq:gevp}
\end{align}
Although one can arguably call \eqref{eq:evp} and \eqref{eq:gevp} 
nonlinear problems they are typically considered linear eigenvalue problems, due
to the linear appearance of $\lambda$ and in
contrast to the quadratic eigenvalue problem \cite{TisK01},
\begin{align*}
  (\lambda^{2}A_{2}+\lambda A_{1} + A_{0})v = 0,
\end{align*}
and the general nonlinear eigenvalue problem \cite{BetHMST13,MehH04},
\begin{align*}
  Q(\lambda)v = 0,
\end{align*}
which both are nonlinear in the eigenvalue. Recently, research interest in
eigenvector nonlinearities, that is, in the basic case
\begin{align*}
  A(v) v = \lambda(v),
\end{align*}
has increased, see, for instance,
\cite{jarlebring2014inverse,cai2018eigenvector,jarlebring2021implicit,claes2021linearizability}.
Problem \eqref{eq:sevp} is different from these. Nonetheless, ideas for the
solution of other (nonlinear) eigenvalue problems are applicable for
\eqref{eq:sevp} as well.

We will discuss two expansions of $A(\mu)$, $x(\mu)$, and $\lambda(\mu)$.  In
Section~\ref{sec:taylor} we will use a truncated Taylor series expansion.  We
start with less powerful Taylor expansion, since this leads to a conceptual
simpler algorithm.  We discuss this algorithm and its complexity in
Subsection~\ref{sec:algorithm:and:complexity}. In
Subsection~\ref{sec:proofofoconcpet} we demonstrate based on an example using
\eqref{eq:example-proof-of-concept} that the algorithm works. We then extend the
approximation to a truncated Chebyshev expansion in
Section~\ref{sec:chebyshev}. Therein we will use a modification of the algorithm
for the Taylor expansion to compute a good starting vector which is then refined
using Newton's method. The paper is concluded with some numerical experiments,
Section~\ref{sec:numerical:experiments}, and conclusions,
Section~\ref{sec:conclusions}.

For our numerical experiments we use Matlab (R2020b) with a machine precision of
$\epsilon_{\text{mach}}\approx 2.2204_{10^{-16}}$ and a computer with Ubuntu
18.04.1, am Intel Core i7-10710U CPU, and 16 GB of RAM.

\section{Taylor expansion}
\label{sec:taylor}
A Taylor expansion of $A(\mu)$ at the expansion point~$\mu_{0}$ is given by
\begin{align*}
  A(\mu) = A_{0} + (\mu-\mu_{0}) A_{1} + \tfrac{1}{2}(\mu-\mu_{0})^{2}A_{2}
  + \tfrac{1}{3!}(\mu - \mu_{0})^{3}A_{3} + \dotsb,
\end{align*}
with
$\left.A_{k}\right|_{ij} = \frac{\partial^{k} A_{ij}}{\partial
  \mu^{k}}(\mu_{0})$. Such an expansion exists, because we assume that
$A(\mu)$ is sufficiently smooth. 
Let us further assume that there are also Taylor expansions for the eigenvector~$v(\mu)$,
\begin{align}
  v(\mu) = v_{0} + (\mu-\mu_{0}) v_{1} + \tfrac{1}{2}(\mu-\mu_{0})^{2}v_{2} + \tfrac{1}{3!}(\mu - \mu_{0})^{3}v_{3} + \dotsb,
  \label{eq:tatylor:v}
\end{align} and for the eigenvalue~$\lambda(\mu)$,
\begin{align*}
  \lambda(\mu) = \lambda_{0} + (\mu-\mu_{0}) \lambda_{1}
  + \tfrac{1}{2}(\mu-\mu_{0})^{2}\lambda_{2}
  + \tfrac{1}{3!}(\mu - \mu_{0})^{3}\lambda_{3} + \dotsb,
\end{align*} with an analogous definition for $v_{k}$ and $\lambda_{k}$.  The
existence of convergent series for $\lambda(\mu)$ and $v(\mu)$ is discussed for
instance in the introduction of Kato's book on perturbation theory for linear
operators \cite{b473}.

We insert these expansions in \eqref{eq:sevp} and compare the coefficients
for $(\mu-\mu_{0})^{p}$. 
Unsurprisingly, the 0th order approximations is equivalent to the solution of
the standard eigenvalue problem obtained at $\mu=\mu_{0}$,
\begin{align}
  A_{0}v_{0} = \lambda_{0}v_{0}.
  \label{eq:0thorder}  
\end{align}
Comparing the next coefficients we find the equation
\begin{align}
  A_{1}v_{0} +A_{0}v_{1} = v_{1} \lambda_{0} + v_{0}\lambda_{1},
  \label{eq:1}
\end{align}
which has to be fulfilled by $v_{1}$ and $\lambda_{1}$.\footnote{Obviously this
  idea is far from novel and can for instance be found in
  \cite{rellich1939storungstheorie}.}  One can reformulate \eqref{eq:1} in
matrix vector form as
\begin{align*}
  \left[\begin{array}{c|c}
    v_{0} & \lambda_{0}I-A_{0}
    \end{array}\right] \begin{bmatrix}
    \lambda_{1}\\v_{1}
  \end{bmatrix} = A_{1}v_{0}.
\end{align*}
This linear system is underdetermined. Although, underdetermined systems can be
solved, for instance with the Moore-Penrose pseudoinverse or Matlab's
\texttt{backslash} operator, we would prefer to avoid an underdetermined
systems. 
We add the condition
\begin{align*}
  1 = v(\mu)^{H}v(\mu),
\end{align*}
which ensures that the eigenvectors are normalized. Using the Taylor series
expansion this condition is
\begin{align*}
  1 = v_{0}^{H}v_{0} + \sum_{k=1}^{p} \frac{1}{k!}\left(\sum_{\ell=0}^{k} \binom{k}{\ell}v_{k-\ell}^{H}v_{\ell}\right)(\mu-\mu_{0})^{k}.
\end{align*}
Thus when solving \eqref{eq:0thorder} we have to normalize the eigenvectors.
For $p=1$ we then have the additional condition
\begin{align*}
  v_{0}^{H}v_{1} = 0.
\end{align*}
This leads to an extension of the linear system to 
\begin{align}
  \label{eq:bordered}
  \left[\begin{array}{c|c}
          0 & v_{0}^{H}\\\hline
          v_{0} & \lambda_{0}I-A_{0}
    \end{array}\right] \begin{bmatrix}
    \lambda_{1}\\v_{1}
  \end{bmatrix} =
  \begin{bmatrix}0\\
    A_{1}v_{0}
  \end{bmatrix},
\end{align} 
which is no longer underdetermined. The equation \eqref{eq:bordered} was used by
Andrew, Chu, and Lancaster to compute the derivative of eigenvectors \cite{ACL93}.

Systems of this form are sometimes referred to as bordered linear systems. If an
efficient solver for $A_{0}$ or $\lambda_{0} I-A_{0}$ is available, for instance
for sparse $A_{0}$, then a block elimination with refinement can be used to solve
\eqref{eq:bordered} efficiently, for details see \cite{govaerts1990block}. 

With $v_{1}$ and $\lambda_{1}$ computed, we can find the linear Taylor
approximation for $\lambda(\mu)$ and $v(\mu)$,
\begin{align*}
  \lambda(\mu) &\approx \lambda_{0} + (\mu-\mu_{0})\lambda_{1},\\
  v(\mu) &\approx v_{0} + (\mu-\mu_{0})v_{1}.\\
\end{align*}
Since we now know $v_{0}$, $v_{1}$, $\lambda_{0}$, and $\lambda_{1}$, we can use
them for the quadratic approximation. Therefore, we have to compare the
coefficients in front of $(\mu-\mu_{0})^{2}$. We have
\begin{align*}
  A_{2}v_{0} +A_{1}v_{1}+A_{0}V_{2} &= v_{2} \lambda_{0}+v_{1}\lambda_{1}+ v_{0}\lambda_{2}
    \intertext{and}
  v_{0}^{H}v_{2} &= -v_{1}^{H}v_{1}.
\end{align*}
Reformulated as a linear system that is 
\begin{align*}
  \left[\begin{array}{c|c}
    0 & v_{0}^{H}\\\hline      
    v_{0} & \lambda_{0}I-A_{0}
    \end{array}\right] \begin{bmatrix}
    \lambda_{2}\\v_{2}
  \end{bmatrix} =
  \begin{bmatrix}
    -v_{1}^{H}v_{1}\\
    A_{2}v_{0} + 2A_{1}v_{1} - 2v_{1}\lambda_{1}
  \end{bmatrix}.
\end{align*}
It is remarkable that this linear system has the same coefficient matrix as before.
Solution of this systems provides $\lambda_{2}$ and $v_{2}$. The approximation is then 
\begin{align*}
  \lambda(\mu) &\approx \lambda_{0} + (\mu-\mu_{0})\lambda_{1}+
                 (\mu-\mu_{0})^{2}\lambda_{2}, \qquad\text{and}\\
  v(\mu) &\approx v_{0} + (\mu-\mu_{0})v_{1} + (\mu-\mu_{0})^{2}v_{2}.
\end{align*}
This can naturally be continued for higher-order Taylor polynomials. We leave it
to the reader to show that
\begin{align}
  \left[\begin{array}{c|c}
          0 & v_{0}^{H}\\\hline
    v_{0} & \lambda_{0}I-A_{0}
    \end{array}\right] \begin{bmatrix}
    \lambda_{k}\\v_{k}
  \end{bmatrix} =
  \begin{bmatrix}-\frac{1}{2}\sum_{\ell=1}^{k-1} \binom{k}{\ell}v_{k-\ell}^{T}v_{\ell}\\
    \sum_{\ell=0}^{k-1} \binom{k}{\ell}A_{k-\ell}v_{\ell} - \sum_{\ell=1}^{k-1}
    \binom{k}{\ell} v_{k-\ell}\lambda_{\ell}
  \end{bmatrix}
  \label{eq:extended:system}
\end{align}
is the linear system determining $\lambda_{k}$ and $v_{k}$. Let $E$ be defined by
\begin{align*}
  E := 
  \left[\begin{array}{c|c}
          0 & v_{0}^{H}\\\hline
    v_{0} & \lambda_{0}I-A_{0}
    \end{array}\right].
\end{align*}

The formulation in \eqref{eq:extended:system} has multiple advantages. We first
discuss the case where one or a few eigenpairs of $A(\mu)$ are computed. If the
matrix $A_{0}$ is symmetric or Hermitian, then so is the coefficient matrix
$E$. If the matrix $A_{0}$ is skew-symmetric/skew-Hermitian, then the sign in
the first row of the linear system should be flipped to preserve the
structure. Since only the diagonal entries of $A_{0}$ are modified to obtain $E$
the zero pattern of most sparse matrices would not be affected
significantly. Thus if $A_{0}$ is a sparse matrix, $E$ will be, too. At most
there are two more nonzeros entries per column in $E$ than there are in
$A$. Sparse direct solvers and sparse iterative solvers can then be
employed. The coefficient matrix $E$, furthermore, is the same for all~$k$. Thus
an $LU$ or Cholesky factorization of $E$ can be reused for all~$k$. Similarly a
preconditioner for $E$ can be reused. If a Krylov subspace method is employed to
solve the linear systems, then it may be possible to recycle selected subspaces
generated for the previous $k$, for details see \cite{Parks2006-fw}. Recycling
subspaces may even be possible between different eigenvalues, i.e.\ for bordered
systems with different shifts in the $\lambda_{0}I-A$ block, if the ideas of
\cite{Soodhalter2016-vs} are extended to bordered systems.

When most or all eigenpairs of $A(\mu)$ are to be computed, one can use an even
more efficient approach. Initially, we need all eigenvalues and eigenvectors of
$A_{0}$. These can all at once be obtained by computing the Schur form of
$A_{0}=QTQ^{H}$ with Francis's implicitly shifted QR
algorithm\footnote{Additional computations are necessary to obtain the
  eigenvectors. Typically this is achieved by inverse iteration or by swapping
  the diagonal entries of $T$. Both require $O(n^{3})$.}  or by computing an
eigenvalue decomposition \mbox{$A_{0}=QDQ^{H}$} with one of the special solvers for
symmetric or skew-symmetric $A_{0}$; both require in general $O(n^{3})$. The
matrix $Q$ is unitary, $T$ upper triangular, and $D$ diagonal. We can now
simplify solving with $E$ by factorizing $E$ with the help of $Q$:
\begin{align}
  E = \left[\begin{array}{c|c}
              0 & v_{0}^{H}\\\hline
              v_{0} & \lambda_{0}I-A_{0}\\
            \end{array}\right] =
  \left[\begin{array}{c|c}
          1 &0\\\hline
          0 & Q\\
        \end{array}\right]
  \underbrace{\left[\begin{array}{c|c}
          0 &v_{0}^{H}Q\\\hline
          Q^{H}v_{0} & \lambda_{0}I-T\\
        \end{array}\right]}_{=\widehat{E}}
  \left[\begin{array}{c|c}
          1 &0\\\hline
          0 & Q^{H}\\
        \end{array}\right].
  \label{eq:decomposition_of_E}
\end{align}
The matrix $\widehat{E}$ is a permutation of an upper-triangular matrix or a
permutation of a diagonal matrix, if $A_{0}$ is symmetric. This turns solving
with $E$ into applying $Q^{H}$ to part of the right-hand side followed by a
backward solve with a permuted upper triangular or diagonal matrix, followed by
another application of $Q$. These steps can be done in $O(n^{2})$ flops.

Unfortunately, there is a big disadvantage as well. The computation of $v_{k}$
and $\lambda_{k}$ depends on all previously computed $v_{\ell}$ and
$\lambda_{\ell}$. Thus one cannot solve the linear system with multiple right
hand sides at the same time, but has to compute them sequentially. Accordingly,
computational errors from earlier steps affect all the subsequent ones. This can
lead to an accumulation of errors in the computed components.

\begin{remark}
  This approach works well for the standard parametric eigenvalue problem
  \eqref{eq:sevp}, since only two expansions are involved on either side of the
  equation. Things get significantly more involved when trying to solve the
  generalized parametric eigenvalue problem \eqref{eq:sgevp}, where three
  expansions are needed on the right-hand side. Thus, we restrict the discussion
  here to the standard parametric eigenvalue problem.
\end{remark}
  
\begin{remark}
  As discussed earlier with
  $v(\mu) = v_{0} + (\mu-\mu_{0})v_{1} + \tfrac{1}{2}(\mu-\mu_{0})^{2}v_{2} +
  \dotsb$ also
  $\gamma v(\mu) = \gamma v_{0} + (\mu-\mu_{0})\gamma v_{1} +
  \tfrac{1}{2}(\mu-\mu_{0})^{2}\gamma v_{2} + \dotsb$, with $\abs{\gamma}=1$, is
  an eigenvector of~\eqref{eq:sevp}. When solving for the 0th order
  approximation one can choose $\gamma$, since $(v_{0},\lambda_{0})$ and
  $(\gamma v_{0},\lambda_{0})$ are both solutions of \eqref{eq:0thorder}.  When
  solving \eqref{eq:bordered} and the subsequent steps the previous choice of
  $\gamma$ forces the solution to be $\gamma v_{1}$ and $\gamma v_{i}$, since
  \begin{align*}
  \left[\begin{array}{c|c}
          0 & \gamma v_{0}^{H}\\\hline
          \gamma v_{0} & \lambda_{0}I-A_{0}
    \end{array}\right] \begin{bmatrix}
    \lambda_{1}\\\gamma v_{1}
  \end{bmatrix} =
  \begin{bmatrix}0\\
    A_{1}\gamma v_{0}
  \end{bmatrix}
  \end{align*}
  is merely a scaling of \eqref{eq:bordered}.  
\end{remark}

We conclude this paragraph by providing a different interpretation useful for the Chebyshev
expansion discussed in Section~\ref{sec:chebyshev}. The comparison of the
coefficients can be represented by the following non-linear block lower
triangular system
\begin{align}
  \label{eq:taylor:nonlin}
  \left[
  \arraycolsep=2.5pt
  \begin{array}{cc|cc|cc|cc}
    0 & v_{0}^{T}&&&&&&\\
    v_{0}& -A_{0}&&&&&&\\\midrule
    0 & v_{1}^{T}& 0 & v_{0}^{T}&&&&\\
    v_{1}& -A_{1}&v_{0}&-A_{0}&&&&\\\midrule
    0 & v_{2}^{T} & 0 & v_{1}^{T}& 0 & v_{0}^{T}&&\\
    \frac{1}{2}v_{2}& -\frac{1}{2}A_{2}&v_{1}& -A_{1}&\frac{1}{2}v_{0}& -\frac{1}{2}A_{0}&&\\\midrule
    0 & v_{3}^{T}& 0 & v_{2}^{T} & 0 & v_{1}^{T}& 0 & v_{0}^{T}\\
    \frac{1}{3!}v_{3}&
    -\frac{1}{3!}A_{3}&\frac{1}{2}v_{2}&-\frac{1}{2}A_{2}&
    \frac{1}{2}v_{1}&-\frac{1}{2}A_{1}&
    \frac{1}{3!}v_{0}&
    -\frac{1}{3!}A_{0}
  \end{array}\right]\left[\begin{array}{c}
    \lambda_{0}\\
    v_{0}\\\midrule
    \lambda_{1}\\
    v_{1}\\\midrule
    \lambda_{2}\\
    v_{2}\\\midrule
    \lambda_{3}\\
    v_{3}    
  \end{array}\right] = \left[\begin{array}{c}
    1\\
    0\\\midrule
    0\\
    0\\\midrule
    0\\
    0\\\midrule
    0\\
    0    
  \end{array}\right],
\end{align}
where $\lambda_{i}$ and $x_{i}$ are the unknowns.
The iterative method derived above solves a diagonally scaled version of this
system by forward substitution. The forward substitution removes the
non-linearity.

\subsection{Algorithm and complexity}
\label{sec:algorithm:and:complexity}
\begin{algorithm2e}[tb] 
  \caption{Approximation of the Parametric Eigenvalue Problem by Taylor Polynomials.}
  \label{algo:taylor1}%
  \KwIn{$A_{0}, A_{1}, A_{2} \dotsc, A_{p}$ with
    $A(\mu) \approx A_{0} + (\mu-\mu_{0})A_{1} +
    \frac{1}{2!}(\mu-\mu_{0})^{2}A_{2} + \dotsb$.} %
  \KwOut{Taylor coefficients for approximations of $\lambda(\mu)$ and $v(\mu)$
    in the neighborhood around $\mu_{0}$.} %
  
  Find one eigenpair $(v_{0},\lambda_{0})$ of $A_{0}$ with
  $A_{0}v_{0}=v_{0}\lambda_{0}$.\\%
  Compute $E =
  \begin{bmatrix}
    0 & v_{0}^{H}\\
    v_{0} & \lambda_{0}I-A_{0}
  \end{bmatrix}$ and, if appropriate, a decomposition of $E$.\\%
  \For{$k=1,\dotsc,p$}{%
    $y := 0$\;%
    $z := 0$\;%
    \For{$\ell=0,\dotsc,k-1$}{%
      $y := y + \binom{k}{\ell} A_{k-\ell}v_{\ell}$\;%
      \If{$\ell>=1$}{%
        $y := y - \binom{k}{\ell} v_{k-\ell}\lambda_{\ell}$\;%
        \If{$\ell<k-1$}{%
          $z = z + \binom{k-1}{\ell-1} v_{k-\ell}^{H}v_{\ell}$\;%
        }%
      }%
    }%
    Solve $E\begin{bmatrix}%
      \lambda_{k}\\v_{k}%
    \end{bmatrix} :=
    \begin{bmatrix}
      -z/2\\y
    \end{bmatrix}$\;%
  }    
\end{algorithm2e} 

Algorithm~\ref{algo:taylor1} shows the steps needed to compute a Taylor
approximation for the eigenpair $(v(\mu),\lambda(\mu))$. This has to be
repeated for each eigenpair of interest. Thus up to $n$ times if all eigenpairs
are required.

We will now assume that $A_{k}$ are dense matrices, and that a matrix-vector
product costs approximately $O(n^{2})$ flops and a matrix-matrix product
$O(n^{3})$ flops. For small to medium size $n$ modern computers with modern
(blocked) implementation typically use algorithms with these computational
complexities.\footnote{Only for (very) large $n$ Strassen-type algorithms are
  occasionally used for matrix-matrix multiplication.}  However, the runtime of
these operations is often limited by the latency of the memory access and the
size of the cache. Thus the runtime of the matrix-matrix multiplication grows
quadratically for many examples with \mbox{$n<100$} or even \mbox{$n<1000$}.

The inner loop of Algorithm~\ref{algo:taylor1} consists of $k$ matrix-vector
products. Computations comparable to an additional matrix-vector product are
required for the solution of the linear system if a decomposition of $E$ is
available. Thus the number of flops for the outer loop can be estimated by
\begin{align*}
  O(p^{2}n^{2}).
\end{align*}
The most expensive part outside the loop is the decomposition of $E$ at
$O(n^{3})$ and the computation of the eigenpair $(v_{0},\lambda_{0})$ also in
$O(n^{3})$ if inverse iteration or Francis's implicitly shifted QR algorithm is
used.\footnote {If Francis's algorithm is used, then all eigenpairs are computed
  in $O(25n^{3})$ \cite{GolV13}.} Only for the largest eigenvalue(s) the power
method at costs of about $O(n^{2})$ can be used.

In total we need $O(n^{3}+p^{2}n^{2})$ flops to compute the Taylor approximation
of degree $p$ of a single eigenpair. When computing all eigenpairs we can make
use of the decomposition \eqref{eq:decomposition_of_E}. That means we need to
compute the Schur decomposition at $O(25n^{3})$ once and $O(p^{2}n^{2})$ flops
per eigenvalue for a total of $O((25+p^{2})n^{3})$ flops.  In the next
subsection we will demonstrate that this algorithm works and that the numerical
experiments do not contradict the complexity estimates.

\subsection{Proof-of-Concept}
\label{sec:proofofoconcpet}

In this section we will discuss a proof-of-concept for Algorithm~\ref{algo:taylor1}.

\begin{example}\label{example:1}
  We pick $n$ points $p_{i}=(x_{i},y_{i},z_{i})$ on a torus with
  \begin{align*}
    x_{i} &= \cos(2\pi\theta_{i})(5 + \cos(4\pi\theta_{i}))\\
    y_{i} &= \sin(2\pi\theta_{i})(5 + \cos(4\pi\theta_{i}))\\
    z_{i} &= \sin(4\pi\theta_{i})
  \end{align*}
  and $\theta_{i}=\frac{i}{n}$, $i=1,\dotsc,n$. These points are aligned on a
  line coiled twice around a torus.  We then define a matrix \textit{$U$} by
  $U_{ij}=\norm{u_{i}-u_{j}}_{2}$ and
  \begin{align*}
    A(\mu) = \exp_{\circ}(-\mu U).
  \end{align*}

  We look at the eigenvalues of $A(\mu)$ with $\mu$ in the interval $[0,1.5]$.
  The eigenvalues intersect in this interval and we are interested in verifying
  that Algorithm~\ref{algo:taylor1} can deal with that.

  At first we pick $n=8$.
\end{example}

We choose $\mu_{0}=0.2$ and compute the 6th order Taylor polynomial to
approximate~$\lambda(\mu)$. The resulting Taylor polynomials are depicted in
Figure~\ref{fig:example1_1_8_7_0.20_1} by the blue lines. The parameter~$\mu$ is plotted
along the $x$-axis, while the non-negative eigenvalues are shown along the
$y$-axis.

\input{code/fig41} %

For comparison we discretize the interval with 151 discretization points
$\hat{\mu}_{i}$ and solve the resulting standard eigenvalue problems for all
$\hat{\mu}_{i}$. The resulting eigenvalues are shown as red crosses. Our goal is
that there is one blue line for each sequence of red crosses ideally going
through the crosses. 

One can observe that near the expansion point $0.2$ the blue lines of the
Taylor approximation follow the red crosses. In particular the behavior at
$0$---one eigenvalue converges to $n=8$, the other seven to $0$---is
represented well.

For $\mu>0.4$ the blue lines show a visible difference from the red crosses. Some
of the blues lines even go below $0$ or above $8$. It is known that for
$\mu\rightarrow\infty$ all eigenvalues converge to $1$. The red crosses
naturally exhibit this convergence while the blue lines do not. This is an
inherent limitation of using polynomials to approximate these functions.

Figure~\ref{fig:example1_1_8_7_0.20_1} also displays a magnified plot around the intersection
of the second and third largest eigenvalue. In the magnified plot the crosses
and the Taylor approximation exhibit the same qualitative behavior. These
intersection points are relevant since the dominant eigenspace change
significantly around them.

The sum of all eigenvalues for a fixed $\mu$ is equal to $n$; $n=8$ in
our example. This can be shown by verifying that the trace is equal to $n$, since
$A_{ii}=0$ for all $i$. As shown in Table~\ref{tab:example1:coefficients} the
Taylor approximations show this behavior, too. Note that we have not used any
additional restrictions on the coefficients enforcing the sum of eigenvalues to
be $8$. One can observe that the other Taylor coefficients add up to
approximately $0$. However, with increasing $p$ the sum of the $\lambda_{p}$
seems to get further and further away from $0$. For $p=10$ the sum is already
$3.78_{10^{-6}}$ and for $p=15$ the sum is $5.25_{10^{+2}}$.  It appears that
the accumulation of errors imposes a limit on the maximum degree for the Taylor
expansion.

We also observe in Table~\ref{tab:example1:coefficients} that the Taylor
coefficients are growing with increasing order.
  
\begin{table}[t]%
  \caption{Coefficients of Taylor series approximations}%
  \centering%
    \include{code/tab41}
  \label{tab:example1:coefficients}
\end{table}
We also timed the algorithm, see Table~\ref{tab:timings}, to check if the
experimental data is in accordance with the complexity, $O(25n^{3}+p^{2}n^{3})$,
discussed in Section~\ref{sec:algorithm:and:complexity}.  We observe timings
in accordance with a complexity of $n^{3}$. In the first two columns we double
$n$ from row to row, but only sometimes observe that the runtime increases
more than fourfold. In the last column we double $p$ from row to row. The
runtime seems to be growing slower than $p^{2}$.

\begin{table}[t]
  \caption{Runtime $t_{i}$ of Algorithm~\ref{algo:taylor1} for the computation of all
    eigenpairs for different combinations of $n$ and $p$.}
  \centering
  \begin{tabular}{rrrr|rrrr|rrrr}
    \toprule
    $n$ & $p$ & $t_{i}$ in s & $\tfrac{t_{i}}{t_{i-1}}$&
    $n$ & $p$ & $t_{i}$ in s & $\tfrac{t_{i}}{t_{i-1}}$&
    $n$ & $p$ & $t_{i}$ in s& $\tfrac{t_{i}}{t_{i-1}}$\\
    \midrule
8 & 2 &   0.0017 & ---&
8 & 7 &   0.0017 & ---&
8 & 2 &   0.0009 & ---\\
16 & 2 &   0.0018 &   1.04 &
16 & 7 &   0.0046 &   2.76 &
8 & 4 &   0.0016 &   1.79 \\
32 & 2 &   0.0041 &   2.30 &
32 & 7 &   0.0078 &   1.70 &
8 & 8 &   0.0038 &   2.37 \\
64 & 2 &   0.0064 &   1.58 &
64 & 7 &   0.0200 &   2.56 &
8 & 16 &   0.0065 &   1.69 \\
128 & 2 &   0.0317 &   4.93 &
128 & 7 &   0.1723 &   8.59 &
8 & 32 &   0.0162 &   2.52 \\
256 & 2 &   0.1349 &   4.25 &
256 & 7 &   0.9903 &   5.75 &
8 & 64 &   0.0595 &   3.66 \\
    \bottomrule\\[0.2ex]
  \end{tabular}
  \label{tab:timings}
\end{table}

\input{code/fig52} 
\input{code/fig53}

In Figure~\ref{fig:example4_1_8_26_0.20_1} we compare the difference between the Taylor
approximation and the eigenvalues computed after discretizing the parametric
eigenvalue problem. We plot the maximum overall eight eigenvalues to declutter
the plot. We observe that the 20th Taylor approximation has an error of about
$10\epsilon_{\text{mach}}$ in the interval $[0.1,0.3]$.

Due to the factor $(\mu-\mu_{0})^{p}$ the higher Taylor coefficients play a
secondary role for the approximation of the eigenvalues for $\mu$ close to the
expansion point. However, one can clearly see that for $\mu=0.6$ and for larger
$\mu$ an increase in the degree does not provide an improvement of the
approximation. This is caused by the accumulation of errors as described above.

To check this hypothesis we have redone the computations for
Figure~\ref{fig:example4_1_8_26_0.20_1} with the only change that the matrix $E$ is rounded to
single precision instead of double precision. The results are shown in
Figure~\ref{fig:example4vpa_1_8_18_0.20_1}. We can observe that higher order approximations do
not improve the quality of the approximation due to the accumulation of errors.

\section{Chebyshev Expansion}
\label{sec:chebyshev}
We observed in the last section that the quality of the Taylor approximations
drops quickly if we get further away from the expansion point. This is expected
but unsatisfactory. Increasing the degree of the Taylor expansion can increase
the interval for which we obtain a satisfactory approximation as seen in
Figure~\ref{fig:example4_1_8_26_0.20_1} and~\ref{fig:example4vpa_1_8_18_0.20_1}.
However, as the preliminary numerical experiments above and the ones in
Section~\ref{sec:numerical:experiments} show, this is limited by the
accumulation of errors. Hence, we need a different approach if we are interested
in an accurate approximation of $\lambda(\mu)$ over a large interval.
We decided to use the Chebyshev expansion of~$A(\mu)$, $\lambda(\mu)$, and~$v(\mu)$.

Let $\{U_{0}(\mu),U_{1}(\mu),U_{2}(\mu),\dotsc \}$ be an orthonormal basis of
Chebyshev polynomials of second kind scaled to the interval
$[\mu_{1},\mu_{2}]$.\footnote{This is an arbitrary choice. Alternatively,
  Chebyshev polynomials of first kind can be used with minor adaptions.}  We can
then write
\begin{align}
  A(\mu) &= A_{0}U_{0}(\mu) + A_{1}U_{1}(\mu) + A_{2}U_{2}(\mu) + \dotsc, \label{eq:def:A}\\
  \lambda(\mu) &= \lambda_{0}U_{0}(\mu) + \lambda_{1}U_{1}(\mu) + \lambda_{2}U_{2}(\mu) + \dotsc,\quad\text{and}\nonumber\\
  v(\mu) &= v_{0}U_{0}(\mu) + v_{1}U_{1}(\mu) + v_{2}U_{2}(\mu) + \dotsc.\nonumber
\end{align}

In Section~\ref{sec:taylor} we assumed that the $A_{i}$ are given. This makes
sense for Example~\ref{example:1} and for other examples since the entry-wise
derivative of $A(\mu)$ can be obtained symbolically in an efficient way. For the
Chebyshev expansion it is an unlikely situation that the $A_{i}$ are given,
i.e.\ that the matrix $A(\mu)$ is given as polynomial in a Chebyshev basis. Thus
we first have to compute them by
\begin{align}
  A_{i} := \langle A(\mu), U_{i}(\mu) \rangle_{U} := \int_{\mu_{1}}^{\mu_{2}} A(\mu)\, U_{i}(\mu)\,
  \frac{2}{\mu_{2}-\mu_{1}}\sqrt{1-\left(\frac{2\mu - \mu_{2}-\mu_{1}}{\mu_{2}-\mu_{1}}\right)^{2}}\dmu.
  \label{eq:sp:cheb}
\end{align}
This can be achieved by a quadrature formula or other standard algorithms. We
choose to use the Chebfun package for these computation in our Matlab
code. In~\eqref{eq:sp:cheb} lies a main difference to the Taylor expansion
approach: When using the Taylor expansion $A_{0}$ is equal to $A(\mu_{0})$,
while when using Chebyshev $A_{0}$ is a weighted average of $A(\mu)$ in
$[\mu_{1},\mu_{2}]$.

Equation \eqref{eq:sp:cheb} defines an inner product in which the Chebyshev
basis of second kind is orthonormal, that is
\begin{align*}
  \langle U_{i}(\mu), U_{j}(\mu) \rangle_{U}=\delta_{ij}.
\end{align*}
We can estimate the approximation error by truncating \eqref{eq:def:A} with the
help of this inner product
\begin{align*}
  A(\mu) =& A_{0}U_{0}(\mu) + A_{1}U_{1}(\mu) + A_{2}U_{2}(\mu) + \dotsc +
           A_{p}U_{p}(\mu) + \Delta_{p}(\mu)\\
  \norm{A(\mu)}_{U}^{2} :=& \langle A(\mu), A(\mu)\rangle_{U} =
  \sum_{i=0}^{p}\langle A(\mu), U_{i}(\mu)\rangle_{U} + \langle A(\mu),\Delta_{p}(\mu)\rangle_{U}  
\end{align*}

The Chebyshev basis is degree graded with the degree of $U_{i}$ being
$i$. Thus the degree of the product $U_{i}U_{j}$ is $i\cdot j$. However, the
product is not equal to $U_{ij}$ as it is the case for the basis
$T_{i}=(\mu-\mu_{0})^{i}$ used for the Taylor expansion. In fact for Chebyshev
polynomials of second kind we have
\begin{align}
  \label{eq:uiuj}
  U_{i}U_{j} = U_{i+j} + U_{i+j-2} + \dotsb + U_{\abs{i-j}+2} + U_{\abs{i-j}}.
\end{align}

If we ignore all but the first term in \eqref{eq:uiuj}, then we obtain an
equation very similar to \eqref{eq:taylor:nonlin}:
\begin{align*}\left[
  \begin{array}{cc|cc|cc|cc}
    0 & v_{0}^{T}&&&&&&\\
    v_{0}& -A_{0}&&&&&&\\\midrule
    0 & v_{1}^{T}& 0 & v_{0}^{T}&&&&\\
    v_{1}& -A_{1}&v_{0}&-A_{0}&&&&\\\midrule
    0 & v_{2}^{T} & 0 & v_{1}^{T}& 0 & v_{0}^{T}&&\\
    v_{2}& -A_{2}&v_{1}& -A_{1}&v_{0}& -A_{0}&&\\\midrule
    0 & v_{3}^{T}& 0 & v_{2}^{T} & 0 & v_{1}^{T}& 0 & v_{0}^{T}\\
    v_{3}&
    -A_{3}&v_{2}&-A_{2}&
    v_{1}&-A_{1}&
    v_{0}&
    A_{0}
  \end{array}\right]\left[\begin{array}{c}
    \lambda_{0}\\
    v_{0}\\\midrule
    \lambda_{1}\\
    v_{1}\\\midrule
    \lambda_{2}\\
    v_{2}\\\midrule
    \lambda_{3}\\
    v_{3}    
  \end{array}\right] = \left[\begin{array}{c}
    1\\
    0\\\midrule
    0\\
    0\\\midrule
    0\\
    0\\\midrule
    0\\
    0    
  \end{array}\right].
\end{align*}
Solving this equation does \emph{not} give us the correct solution because we
made a severe simplification.  However, this is an approximation to the
solution, which can be computed iteratively and similar to
Algorithm~\ref{algo:taylor1} due to the block lower triangular structure. We can
obtain this approximation also in $O(25n^{3} + p^{2}n^{2})$ for one eigenpair or
$O((25+p^{2})n^{3})$ for all eigenpairs.

If we do not ignore the terms in~\eqref{eq:uiuj}, then the system of equations
for $p=3$ is
\begin{align}\left[
  \begin{array}{cc|cc|cc|cc}
    0 & v_{0}^{T}&0&{\color{SPECred}v_{1}^{T}}&&&&\\
    v_{0}& -A_{0}&{\color{SPECred}v_{1}}&{\color{SPECred}-A_{1}}&&&&\\\midrule
    0 & v_{1}^{T}& 0 & v_{0}^{T}+{\color{SPECblue}v_{2}^{T}}&0&{\color{SPECorange}v_{1}^{T}}&&\\
    v_{1}& -A_{1}&v_{0}+{\color{SPECblue}v_{2}}&-A_{0}{\color{SPECblue}-A_{2}}&{\color{SPECorange}v_{1}}&{\color{SPECorange}-A_{1}}&&\\\midrule
    0 & v_{2}^{T} & 0 & {\color{SPECred}v_{1}^{T}}& 0 & v_{0}^{T}&&\\
    v_{2}& -A_{2}&{\color{SPECred}v_{1}}& {\color{SPECred}-A_{1}}&v_{0}& -A_{0}&&\\\midrule
    0 & v_{3}^{T}& 0 & {\color{SPECblue}v_{2}^{T}} & 0 & {\color{SPECorange}v_{1}^{T}}& 0 & v_{0}^{T}\\
    v_{3}&-A_{3}&{\color{SPECblue}v_{2}}&{\color{SPECblue}-A_{2}}& {\color{SPECorange}v_{1}}&{\color{SPECorange}-A_{1}}& v_{0}& A_{0}
  \end{array}\right]\left[\begin{array}{c}
    \lambda_{0}\\
    v_{0}\\\midrule
    \lambda_{1}\\
    v_{1}\\\midrule
    \lambda_{2}\\
    v_{2}\\\midrule
    \lambda_{3}\\
    v_{3}    
  \end{array}\right] = \left[\begin{array}{c}
    1\\
    0\\\midrule
    0\\
    0\\\midrule
    0\\
    0\\\midrule
    0\\
    0    
  \end{array}\right], \label{eq:nonlin:sy}
\end{align}
where the colored blocks occur multiple times. The system \eqref{eq:nonlin:sy}
is no longer block lower triangular and thus cannot be solved by forward
substitution. We hence fallback on solving this non-linear system by Newton's
method. In each step we have to compute and invert the Jacobi matrix. This is a
matrix of size $(p+1)(n+1)$. Thus the costs are in $O(n^{3}p^{3})$ for each
Newton step for each eigenvalue. We observe that 4 to 8 Newton steps are
typically sufficient, since we already start with a good approximation.

We truncate the series for $A(\mu)$, $\lambda(\mu)$, and $v(\mu)$ all after
$p+1$ terms, so that they are all polynomials of degree $p$. Different choice
may be possible but we did not see an advantage in using different degrees for
$A$, $\lambda$, and $v$ in our preliminary experiments.

Nonlinear systems can have many solutions. However, the special constructions of
\eqref{eq:nonlin:sy} ensures that every of its solutions represents a Chebyshev
approximation of an eigenpair of the parametric eigenvalue problem. For the
solution
\begin{align*}
  &\begin{bmatrix}
    \lambda_{0} & \phantom{\gamma}v_{0} & \lambda_{1} & \phantom{\gamma}v_{1} & \dotsb & \lambda_{p-1} & \phantom{\gamma}v_{p-1}
  \end{bmatrix}^{T}
\intertext{%
  there are infinitely many more of the form:%
}
 &\begin{bmatrix}
  \lambda_{0} & \gamma v_{0} & \lambda_{1} &\gamma v_{1} & \dotsb & \lambda_{p-1} & \gamma v_{p-1}
\end{bmatrix}^{T} 
\end{align*}
for all $ \gamma$ with $\abs{\gamma}=1$. Newton's method does not guarantee that
we converge to the closest root of these. Thus, despite starting with $n$ initial
approximations, one close to each root, there is no guarantee that we end up
with $n$ distinct eigenpairs $(\lambda(\mu),\gamma v(\mu))$ not just differing
in $\gamma$.

For polynomial rootfinding this problem can be overcome by the Ehrlich-Aberth
method \cite{aberth1973iteration,ehrlich1967modified}, see for instance
\cite{bini2014solving}. This is possible for polynomial rootfinding since the
eigenvectors corresponding to each root can be derived easily from said
root. Unfortunately, this is not possible here. As a consequence we are left
with the hope that the starting vectors are sufficiently close. In the special
case that all eigenvalues and eigenvectors are real, there are only two choices
for $\gamma$, $\gamma=1$ and $\gamma=-1$. In this case we had some heuristic
success in employing the Ehrlich-Aberth method. However, with poorly chosen
start-vectors it was still possible to trick the algorithm into finding the same
eigenpair twice.

\subsection{Accuracy}
We now want to discuss the accuracy of the computed and truncated Chebyshev
expansions for $\lambda(\mu)$ and $v(\mu)$. The Chebyshev theorem states that
the best approximation in the Chebyshev basis $b_{n}$ fulfills
\begin{align}
  f(x) - b_{n}(x) = \frac{f^{(n+1)}(\xi)}{(n+1)!} \prod_{i=0}^{n}(x-z_{i}), \label{eq:cheb:thm}
\end{align}
for unknown interpolation points $z_{i}\in [\mu_{1},\mu_{2}]$ and an unknown
point $\xi\in[\mu_{1},\mu_{2}]$ \cite{Dzy95}. The reason is that the Chebyshev
approximation $b_{n}(x)$ intersects with $f(x)$ at $n$ unknown points in
$[\mu_{1},\mu_{2}]$. Thus it is an interpolation polynomial at these points and
for such an interpolation polynomial \eqref{eq:cheb:thm} holds. We can bound the
right-hand side by
$(b-a)^{n} \sup_{\xi\in[\mu_{1},\mu_{2}]} \frac{f^{(n+1)}(\xi)}{(n+1)!}$. This
is a good bound for the eigenvalues, but for eigenvectors the derivative may be
huge when two eigenvalues are near to each other, see \cite[Sect. 7.2.2]{GolV13}
for details.

Additionally, the Chebyshev coefficients we compute with Newton's method are
perturbed because we only approximate the solution of \eqref{eq:nonlin:sy}.

\subsection{Proof-of-Concept}
We repeat the experiments from Section~\ref{sec:proofofoconcpet} based on
Example~\ref{example:1} with Chebyshev expansion instead of Taylor series
expansion.

\input{code/fig52_cheb}

Figure~\ref{fig:example4_0_8_20_0.25_1} is the Chebyshev version of
Figure~\ref{fig:example4_1_8_26_0.20_1} with Chebyshev expansions on the
interval $[\tfrac{1}{4},1]$. Contrary to the Taylor approximation we do not just
see a small error near the expansion point but a small error over the whole
interval. The error increases slightly towards the endpoints of the intervals,
with the error near $1$ being $10^{-13}$, while the error near $\tfrac{1}{4}$ is
closer to $10^{-14}$. Outside the approximation interval the error increases the
further we are from the endpoints. Here, we do not observe that an increase in
the order leads to worse results, since the use of Newton's method to solve the
nonlinear equation gets us around the error accumulation problem. However, in
Figure~\ref{fig:example4_0_8_30_0.50_2} depicting the error in the eigenvectors,
we observe that the error outside the approximation interval increases for large
$p$.

\section{Numerical Experiments}
\label{sec:numerical:experiments} 
In the last section we used Example~\ref{example:1} to verify that the proposed
methods work for the arguably most easy problem of symmetric positive definite
matrices. In this section we will present further experiments for other
examples.

\begin{figure}
  \centering
  \scalebox{0.8}{
  \begin{tikzpicture}

    \fill [pattern = north west lines] (0,-1) rectangle (-0.2,1);
    \draw [thick] (0,-1) -- (0,1);

    \draw [decoration ={%
      coil,%
      segment length = 1mm,%
      amplitude = 2mm,%
      aspect = 0.5,%
      post length = 2mm,%
      pre length = 2mm},%
    decorate] (0,0) -- (1.5,0) node [midway, above =1.0ex]{$1$};
      
    \node[draw, minimum width = 1cm, minimum height = 0.75cm, anchor=west] (m1) at (1.5,0) {1};

    \draw [decoration ={%
      coil,%
      segment length = 1mm,%
      amplitude = 2mm,%
      aspect = 0.5,%
      post length = 2mm,%
      pre length = 2mm},%
    decorate] (2.5,0) -- (4.0,0) node [midway, above =1.0ex]{$1$};
    
    \node[minimum width = 1cm, minimum height = 0.75cm, anchor=west] (m2) at (4.0,0) {$\dotsm$};
    
    \draw [decoration ={%
      coil,%
      segment length = 1mm,%
      amplitude = 2mm,%
      aspect = 0.5,%
      post length = 2mm,%
      pre length = 2mm},%
    decorate] (5.0,0) -- (6.5,0) node [midway, above =1.0ex]{$1$};

    \node[draw, minimum width = 1cm, minimum height = 0.75cm, anchor=west] (m3) at (6.5,0) {$\mu$};

    \draw [decoration ={%
      coil,%
      segment length = 1mm,%
      amplitude = 2mm,%
      aspect = 0.5,%
      post length = 2mm,%
      pre length = 2mm},%
    decorate] (7.5,0) -- (9.0,0) node [midway, above =1.0ex]{$1$};

    \node[draw, minimum width = 1cm, minimum height = 0.75cm, anchor=west] (m4) at (9.0,0) {$\mu$};

    \draw [decoration ={%
      coil,%
      segment length = 1mm,%
      amplitude = 2mm,%
      aspect = 0.5,%
      post length = 2mm,%
      pre length = 2mm},%
    decorate] (10.0,0) -- (11.5,0) node [midway, above =1.0ex]{$1$};

    \node[minimum width = 1cm, minimum height = 0.75cm, anchor=west] (m5) at (11.5,0) {$\dotsc$};

    \draw [decoration ={%
      coil,%
      segment length = 1mm,%
      amplitude = 2mm,%
      aspect = 0.5,%
      post length = 2mm,%
      pre length = 2mm},%
    decorate] (12.5,0) -- (14.0,0) node [midway, above =1.0ex]{$1$};

    \fill [pattern = north east lines] (14,-1) rectangle (14.2,1);
    \draw [thick] (14,-1) -- (14,1);

  \end{tikzpicture}}
  \caption{Example~\ref{example:2}---Sketch of the springs and masses.}%
  \label{fig:example2}%
\end{figure}
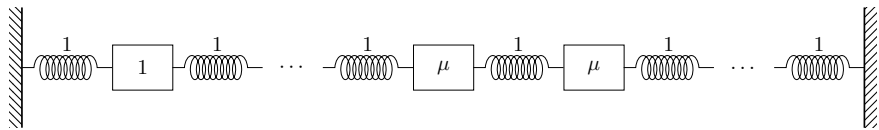
\begin{example}\label{example:2} 
  We investigate a sequence of masses connected with springs as in
  Figure~\ref{fig:example2}.  All springs have a stiffness of~$1$. All but two
  masses in the middle have mass~$1$. The two masses in the middle are both of
  mass~$\mu$.  This example is an extreme simplification of a passenger of
  unknown mass in a car.

  This problem leads to a generalized eigenvalue problem dependent on~$\mu$. As
  mentioned earlier generalized eigenvalue problems are significantly more
  difficult to handle with the expansion approaches. Thus we choose to turn the
  problem into a standard eigenvalue problem. Inversion of the mass matrix leads
  to
  \begin{align*}
    \setlength\arraycolsep{4pt}
    A(\mu) =
    \begin{bmatrix}
      1\\
      &\smash{\ddots}\\
      &&1\\
      &&&\mu\\
      &&&&\mu\\
      &&&&&1\\
      &&&&&&\smash{\ddots}\\
      &&&&&&&1\\      
    \end{bmatrix}^{-1}
    \begin{bmatrix}
      2&-1\\
      -1&\smash{\ddots}&\smash{\ddots}\\
      &\smash{\ddots}&2&-1\\
      &&-1&2&-1\\
      &&&-1&2&-1\\
      &&&&-1&2&\smash{\ddots}\\
      &&&&&\smash{\ddots}&\smash{\ddots}&-1\\
      &&&&&&-1&2\\
    \end{bmatrix}.
  \end{align*}
  Thus only two rows of $A(\mu)$ are dependent on $\mu$.

  The matrix $A(\mu_{0})$ is not symmetric in this example. However, using the
  square root of the mass matrix would permit to construct a symmetric problem
  instead.
\end{example}

We apply Algorithm~\ref{algo:taylor1} with $\mu_{0}=0.8$ and its Chebyshev
sibling for the interval $[\mu_{1},\mu_{2}]=[0.5,1.0]$ to
Example~\ref{example:2}. The resulting eigenvalues are depicted in
Figure~\ref{fig:example1_1_8_7_0.80_2}. 
Figure~\ref{fig:example4_1_8_26_0.80_2}
shows the error of the Taylor approximation to sampled
eigenvalues. Figure~\ref{fig:example4_0_8_30_0.50_2} does the same for the
Chebyshev approximation. We observe that in this example we need a higher order
for the Chebyshev approximation than for the Taylor approximation in order to
obtain the same level of accuracy. However, the Chebyshev approximation is
better over a larger interval. The example demonstrates that the technique
computing the Chebyshev approximation does not suffer from an accumulation of
errors and thus we can compute the higher order approximation sufficiently
accurate. It also shows that the Chebyshev approximation is of good quality in
the inner part of the interval $[\mu_{0},\mu_{1}]$. However, near the endpoints
the quality deteriorates significantly.

\input{code/fig61_both} 
\input{code/fig62} 
\input{code/fig62_cheb}

The examples we looked at so far had a full set of eigenvectors for all
parameter values in the interval of interest. The following examples uses a
basic Jordan block to test the algorithm for an example with defective
eigenvalues. 

\begin{example}
  \label{example:3}
  We use a Jordan block with a parameter in the lower left corner, that is
  \begin{align*}
    A(\mu) =
    \begin{bmatrix}
      1 & 1 & &&0\\
      0 & 1 & 1\\
      \vdots&\ddots & \ddots &\ddots\\
      0&\cdots & 0& 1& 1\\
      \mu &0&\cdots & 0 &1
    \end{bmatrix}
  \end{align*}
  The eigenvalues of $A(\mu)$ are the roots of the characteristic polynomial
  \begin{align*}
    p(\lambda,\mu) = (\lambda-1)^{n}-\mu.
  \end{align*} The roots are the unit roots $\rho_{1},\dotsc,\rho_{n}$ of
  $\rho^{n}=-1$ scaled by $\mu^{1/n}$ and shifted by $1$ to the right.
\end{example}

We note that we primarily focused our algorithms on the case of non-defective,
preferable positive definite matrices, which Example~\ref{example:3} is very
much not. Nevertheless we are interested in how far the techniques described in
this paper produce meaningful results in this case.

We apply Algorithm~\ref{algo:taylor1} with $\mu_{0}=0.2$ and $\mu_{0}=-0.2$ and
its Chebyshev sibling with $[\mu_{1},\mu_{2}]=[0.5,1.0]$ to
Example~\ref{example:3}. The real and imaginary parts of the Taylor
approximation eigenvalues are depicted in
Figure~\ref{fig:example1_1_8_26_0.20_3} for $\mu_{0}=0.2$. A similar figure can
be obtained 
for $\mu_{0}=-0.2$. We observe that the Taylor series approximations do not
provide useful approximations beyond the singularity at $\mu=0$. This is
expected, since the eigenvalues are not analytic at
$\mu=0$. Algorithm~\ref{algo:taylor1} failed for the expansion point
$\mu_{0}=0$.

Using the Chebyshev approximation we observe a very similar behavior, see
Figure~\ref{fig:example1_0_8_20_0.10_3}. The algorithm failed to provide a
meaningful approximation when $0\in[\mu_{1},\mu_{2}]$. The Chebyshev
approximation is also not capable of approximating the eigenvalues beyond $0$.

\input{code/fig71}
\input{code/fig71_cheb}

\subsection{Assessing the Quality of the Eigenvector Approximation}
We are now going to assess the quality of the eigenvector approximation. We use
the Taylor and Chebyshev approximation of the eigenvectors, respectively. We
observe that despite our efforts to produce normalized eigenvectors, the
function values $v(\mu)$ are not of norm $1$. In particular, the
Chebyshev approximation techniques produces eigenvectors far from norm $1$. For
our comparison we evaluate the function $v(\mu)$ for a particular $\mu$ and then
normalize the result before comparing it to the normalized eigenvectors of
$A(\mu)$.

We fix $\mu_{j}$. Let $V_{s}$ be the matrix of eigenvectors of $A(\mu_{j})$ with
$A(\mu_{j})V_{s}=D_{s}V_{s}$, with $D_{s}$ diagonal, and let $V$ be the result
of the normalization of $v(\mu_{j})$. We use the Matlab command
\texttt{max(abs(max(abs(Vs'*V))-1))} to compute the maximum deviation of the
eigenvectors for $\mu_{j}$.

Figure~\ref{fig:example5_1_8_26_0.20_1}--\ref{fig:example5_0_8_30_0.50_2} show
the deviation of the eigenvectors for Example~\ref{example:1} and~\ref{example:2}. We observe in
Figure~\ref{fig:example5_1_8_26_0.20_1} and
Figure~\ref{fig:example5_0_8_30_0.50_2} that increasing the order of
approximation provides diminishing returns and eventually a decrease in
approximation quality.

Surprisingly, Figure~\ref{fig:example5_0_8_20_0.25_1} shows that order 10 is
sufficient to approximate the eigenvectors with the Chebyshev approach well,
while Figure~\ref{fig:example4_0_8_20_0.25_1} shows that we need a higher order
to approximate the eigenvalues well. Knowing the eigenvectors $v(\mu_{i})$ for a
given $\mu_{i}$ and a known matrix $A(\mu_{i})$ allows to compute an
approximation to the eigenvalues $\lambda(\mu_{i})$ using the Rayleigh
quotient. It is a well known fact that for symmetric matrices an approximation
$q$ to the eigenvector $v$ with $\norm{q-v}=O(\delta)$ implies that the Rayleigh
quotient $\rho=q^{T}Aq/q^{T}q$ approximates the eigenvalue $\lambda$ to $v$ with
$\abs{\lambda-\rho}=O(\delta^{2})$, see for instance \cite[Exercise
5.4.33]{b148ed3}. Thus for symmetric matrices using the eigenvector approximation
together with the Rayleigh quotient is very likely going to produce a better
approximation to the eigenvalue. This is aided by the fact that for the Rayleigh
quotient one can use the exact $A(\mu)$. We tried this out and found a better
approximation than using $\lambda(\mu_{i})$, cp.\
Figure~\ref{fig:example7_0_8_20_0.25_1} with
Figure~\ref{fig:example4_0_8_20_0.25_1}.

\input{code/fig61_taylor_1}
\input{code/fig61_taylor_2}
\input{code/fig61_cheb_1}
\input{code/fig61_cheb_2}

\input{code/fig71_cheb_1}

\subsection{Sampling}
We have presented a method to find approximations to $\lambda(\mu)$. These
approximations can be used to speed-up the following sampling process: Assume a
distribution of $\mu$ is given and can be sampled. For the following tests we
use $\mu\sim \cN(\mu_{0},0.01)$, that is normally distributed around
$\mu_{0}=0.2$ with standard deviation $0.1$. We can thus produce samples
$\mu_{i}$ and use these to compute samples $\lambda(\mu_{i})$ of the
distribution of the eigenvalues. In the example below the second and third
largest eigenvalue are sampled. This is faster than computing $A(\mu)$ and then
using a standard eigenvalue algorithm to compute the eigenvalues of
$A(\mu)$. For $10,000$ samples using Example~\ref{example:1} with $n=8$ and
approximations of $6$th order we observed that $3.6938$ s were spend on forming
the matrix as a function, $0.0096$ s were spend on solving \eqref{eq:sevp} with
the Taylor approximation and $0.1464$ s with the Chebyshev approximation. For
sampling the eigenvalues using the full matrix $A(\mu)$ $1.6386$ s were needed,
while sampling with the obtained Taylor (using Horner's method) and Chebyshev
approximations could be done in $0.0807$ s and $0.0616$ s, respectively. Thus
the sampling was done 20.3 and 26.6 times faster, respectively, or 18.14 and 7.87
times faster if the one-time costs for solving \eqref{eq:sevp} are added to the
costs for the sampling. The more accurate results obtained through 
the Rayleigh quotient are slightly more expensive and, thus, are only 5.58 and
5.89 times faster, respectively. We summarize these numbers together with the
theoretical complexities in Table~\ref{tab:sampling}.

\begin{table}[t]
  \caption{Comparing direct sampling of $s=10,000$ eigenvalues of $A(\mu)$ with first
    computing a Taylor or Chebyshev expansion, for $A(\mu)\in\Rnn$, with $n=8$ and
    $p=6$.}
  \centering
  \begin{tabular}{ccccc}
    \toprule
    && Taylor & Chebyshev & direct sampling\\
    \midrule   
    \multirow{2}{*}{\rotatebox{90}{setup}}&theoretical costs & $O((25+p^{2})n^{3})$ & $O(n^{4})$ & ---\\
    &runtime & 0.0096 s & 0.1464 s & ---\\
    \midrule
    \multirow{5}{*}{\rotatebox{90}{sampling}} & theoretical costs & $O(ps)$ & $O(ps)$ &
                                                                          $O(sn^{3})$\\
    &runtime & 0.0807 s & 0.0616 s & 1.6386 s\\
    &speed up & 20.3$\times$ & 26.6$\times$ & ---\\
    &with Rayleigh quot. & 0.2838 s & 0.1319 s & ---\\
    &speed up & 5.77$\times$ & 12.4$\times$ & ---\\
    \midrule
    \multirow{4}{*}{\rotatebox{90}{combined}}                                                                      
    &runtime & 0.0903 s & 0.2080 s & 1.6386 s\\
    &speed up & 18.14$\times$ & 7.87$\times$ & ---\\
    &with Rayleigh quot. & 0.2934 s & 0.2783 s & ---\\
    &speed up & 5.58$\times$ & 5.89$\times$ & ---\\
    
    \bottomrule\\[0.2ex]
  \end{tabular}
  \label{tab:sampling}
\end{table}

\begin{figure}
  \centering
  \includegraphics{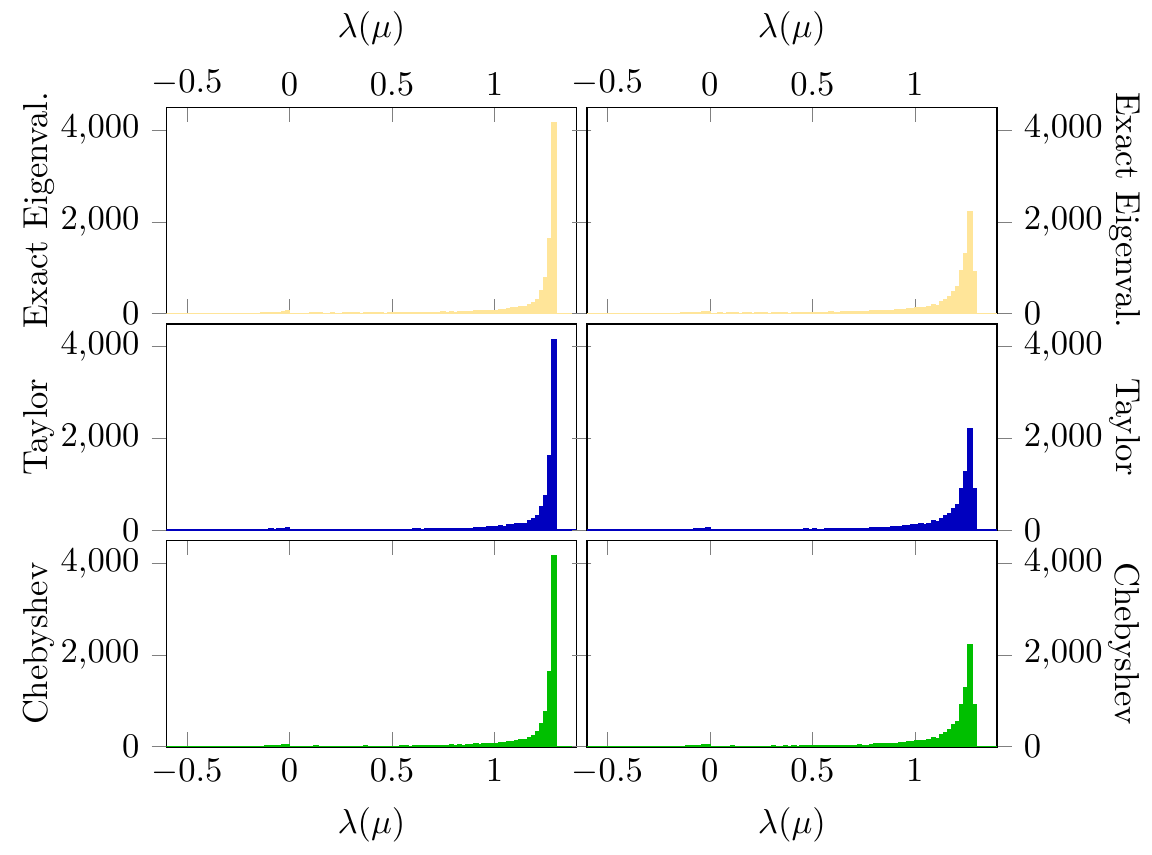}
\caption{Example~\ref{example:1}---Histogram plot of eigenvalue distribution second and third
  largest eigenvalue, $6$th order Taylor approximation with $\mu_{0}=0.20$ and
  $6$th order Chebyshev approximation with $[\mu_{1},\mu_{2}]=[0.1,0.3]$.}%
\label{fig:histogram1}%
\end{figure}

With the help of the sampling we generated the histogram plots in
Figure~\ref{fig:histogram1}. The histogram plots appear visually very similar
and can be used to obtain a general impression of the distribution. In
particular for $\mu$ far from $\mu_{0}$ or from $[\mu_{1},\mu_{2}]$ the
approximation quality is low.

\begin{figure}
  \centering
  \includegraphics{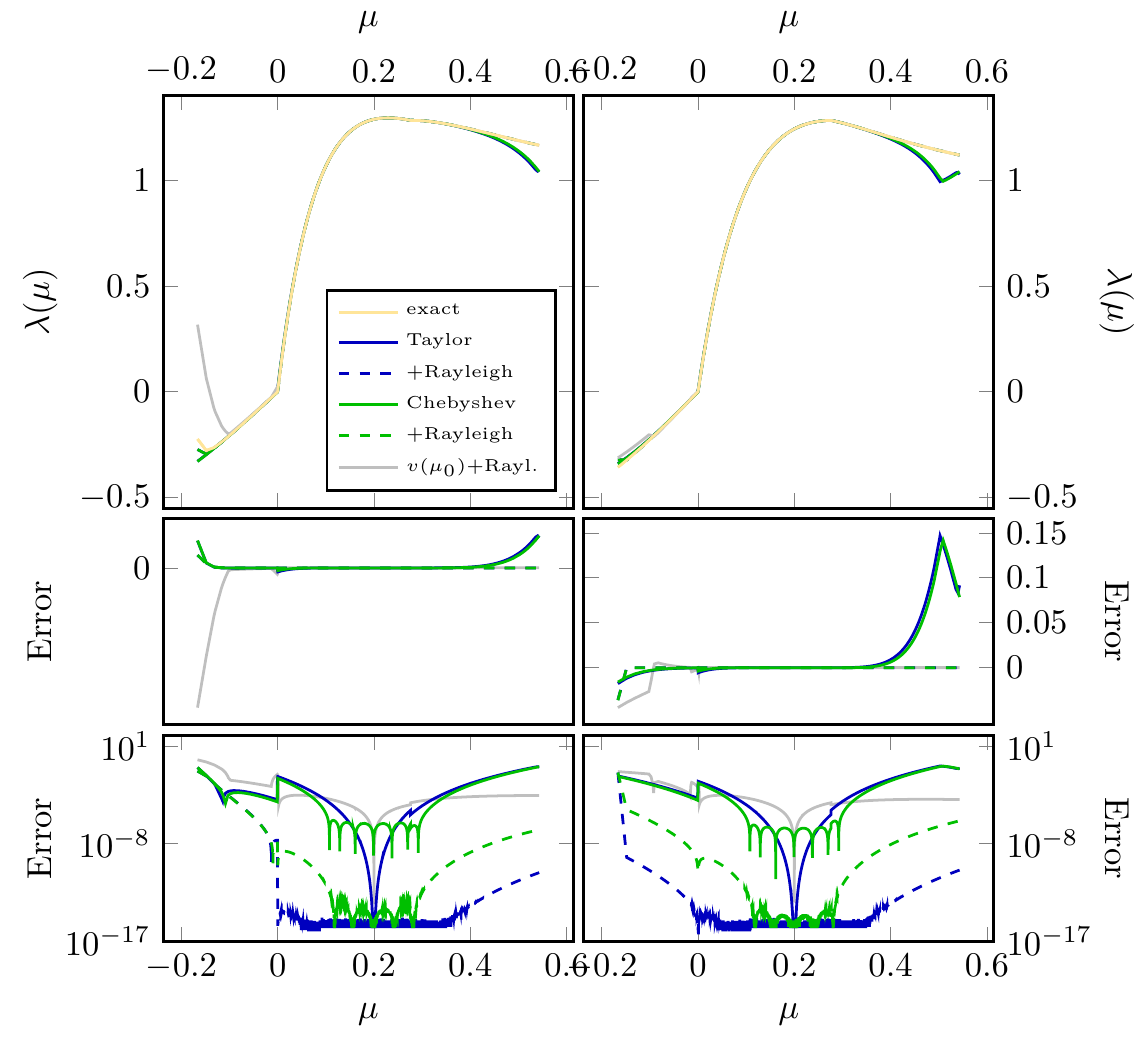}
\caption{Example~\ref{example:1}---$\lambda(\mu)$ over $\mu$ second (left) and
  third largest eigenvalue (right), exact eigenvalue
  (orange line \protect\begin{tikzpicture}[x=1cm,y=1cm]
    \protect\draw[very thick,SPECCorange] (0,0)--(0.4,0); \protect\node[minimum
    height=0,minimum width=0] at (0.2,0) {}; \protect\end{tikzpicture}\,), $6$th
  order Taylor approximation with $\mu_{0}=0.20$ (solid
  blue line \protect\begin{tikzpicture}[x=1cm,y=1cm]
    \protect\draw[very thick,SPECblue] (0,0)--(0.4,0); \protect\node[minimum
    height=0,minimum width=0] at (0.2,0) {}; \protect\end{tikzpicture}\,), and
  $6$th order Chebyshev approximation eigenvalue (solid
  green line \protect\begin{tikzpicture}[x=1cm,y=1cm]
    \protect\draw[very thick,SPECgreen] (0,0)--(0.4,0); \protect\node[minimum
    height=0,minimum width=0] at (0.2,0) {}; \protect\end{tikzpicture}\,) with
  $[\mu_{1},\mu_{2}]=[0.1,0.3]$; and eigenvalue approximation using the Rayleigh
  quotient based on the eigenvector approximation (dashed lines
  \protect\begin{tikzpicture}[x=1cm,y=1cm] \protect\draw[very
    thick,SPECblue,dashed] (0,0)--(0.4,0); \protect\node[minimum
    height=0,minimum width=0] at (0.2,0) {};
    \protect\end{tikzpicture}\,and\,\protect\begin{tikzpicture}[x=1cm,y=1cm]
    \protect\draw[very thick,SPECgreen,dashed] (0,0)--(0.4,0);
    \protect\node[minimum height=0,minimum width=0] at (0.2,0) {};
    \protect\end{tikzpicture}\,). For comparison the Rayleigh quotient
  $v(\mu_{0})^{T}A(\mu)v(\mu_{0})/v(\mu_{0})^{T}v(\mu_{0})$
  (grey line \protect\begin{tikzpicture}[x=1cm,y=1cm]
    \protect\draw[very thick,SPECblack,dashed] (0,0)--(0.4,0);
    \protect\node[minimum height=0,minimum width=0] at (0.2,0) {};
    \protect\end{tikzpicture}\,).}%
 \label{fig:curve1}%
\end{figure}

\begin{figure}
  \centering
  \includegraphics{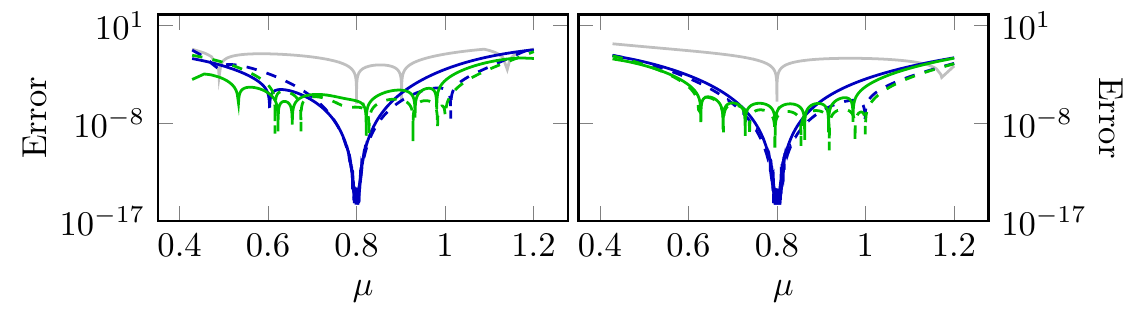}
\caption{Example~\ref{example:2}---$\lambda(\mu)$ over $\mu$ second (left) and
  third largest eigenvalue (right),  $6$th order Taylor approximation with
  $\mu_{0}=0.80$ (solid
  blue line \protect\begin{tikzpicture}[x=1cm,y=1cm]
    \protect\draw[very thick,SPECblue] (0,0)--(0.4,0); \protect\node[minimum
    height=0,minimum width=0] at (0.2,0) {}; \protect\end{tikzpicture}\,), and $6$th order
  Chebyshev approximation eigenvalue (solid
  green line \protect\begin{tikzpicture}[x=1cm,y=1cm]
    \protect\draw[very thick,SPECgreen] (0,0)--(0.4,0); \protect\node[minimum
    height=0,minimum width=0] at (0.2,0) {}; \protect\end{tikzpicture}\,))
  with $[\mu_{1},\mu_{2}]=[0.6,1.0]$; and eigenvalue approximation using the
  Rayleigh quotient based on the eigenvector approximation (dashed lines
  \protect\begin{tikzpicture}[x=1cm,y=1cm] \protect\draw[very
    thick,SPECblue,dashed] (0,0)--(0.4,0); \protect\node[minimum
    height=0,minimum width=0] at (0.2,0) {};
    \protect\end{tikzpicture}\,and\,\protect\begin{tikzpicture}[x=1cm,y=1cm]
    \protect\draw[very thick,SPECgreen,dashed] (0,0)--(0.4,0);
    \protect\node[minimum height=0,minimum width=0] at (0.2,0) {};
    \protect\end{tikzpicture}\,). For
  comparison the Rayleigh quotient
  $v(\mu_{0})^{T}A(\mu)v(\mu_{0})/v(\mu_{0})^{T}v(\mu_{0})$
  (grey line \protect\begin{tikzpicture}[x=1cm,y=1cm]
    \protect\draw[very thick,SPECblack,dashed] (0,0)--(0.4,0);
    \protect\node[minimum height=0,minimum width=0] at (0.2,0) {};
    \protect\end{tikzpicture}\,).}%
\label{fig:curve2}%
\end{figure}

We also plotted $\lambda(\mu)$ over $\mu$ (\,\begin{tikzpicture}[x=1cm,y=1cm]
  \draw[very thick,SPECCorange] (0,0)--(0.4,0);
  \node[minimum height=0,minimum width=0] at (0.2,0) {};
\end{tikzpicture}\,) in Figure~\ref{fig:curve1}
together with the error. We observe that the Taylor approximation (\,\begin{tikzpicture}[x=1cm,y=1cm]
  \draw[very thick,SPECblue] (0,0)--(0.4,0);
  \node[minimum height=0,minimum width=0] at (0.2,0) {};
\end{tikzpicture}\,) provides a
good approximation near $\mu_{0}$, while the Chebyshev approximation (\,\begin{tikzpicture}[x=1cm,y=1cm]
  \draw[very thick,SPECgreen] (0,0)--(0.4,0);
  \node[minimum height=0,minimum width=0] at (0.2,0) {};
\end{tikzpicture}\,) is better over the whole interval $[\mu_{1},\mu_{2}]$. We can again
observe that the Rayleigh quotient can be used to improve the approximation to
$\lambda(\mu)$ for both the Taylor and the Chebyshev approximation; dashed
lines \begin{tikzpicture}[x=1cm,y=1cm]
  \draw[very thick,SPECblue,dashed] (0,0)--(0.4,0); \node[minimum
  height=0,minimum width=0] at (0.2,0) {};
\end{tikzpicture}\,and\,\begin{tikzpicture}[x=1cm,y=1cm]
  \draw[very thick,SPECgreen,dashed] (0,0)--(0.4,0);
  \node[minimum height=0,minimum width=0] at (0.2,0) {};
\end{tikzpicture}, respectively. However, this improvement only works for symmetric matrices,
like Example~\ref{example:1}. Figure~\ref{fig:curve2} presents the last row of
Figure~\ref{fig:curve1} but for the non-symmetric Example~\ref{example:2}, and
hence, there is only a minor or no improvement visible when using the Rayleigh
quotient approximation.

\section{Conclusions}
\label{sec:conclusions}
We have presented a Taylor approximation based method to compute approximations
to $\lambda(\mu)$ and $v(\mu)$ for the parametric eigenvalue problem
\eqref{eq:sevp}. The presented algorithm works for small degrees $p$ and the
investigated examples. The algorithm accumulates the errors when increasing the
degree $p$. Together with the errors present in $E$ when forming the matrix in
double (or single) precision, this leads to a serious limitation regarding the
maximum degree $p_{\max}$. This limits the usefulness of this approach to a
small neighborhood around $\mu_{0}$. We were able to verify that the runtime
complexity of this algorithm is within the expectations set by counting the
number of flops. The algorithm computes all eigenvalues in
$O((25+p^{2})n^{3})$.

To overcome the limitations we extended the approach to Chebyshev
approximation. This requires the solution of a non-linear system with Newton's
method. The Chebyshev approach has higher costs. Despite a good starting point
for the Newton iteration, there is no guarantee that all eigenvalues can be
found. However, in the experiments the Chebyshev approximations are superior to
the Taylor approximations, in particular since a high accuracy can be achieved
over a given interval and not just in the neighborhood of an expansion point. We
showed that the method can be used for the sampling of eigenvalues if a
distribution of the parameter is given. Depending on the number of sampling
points the method presented here can be significantly faster than Monte-Carlo
methods.

\section*{Acknowledgments}\label{sec:ack}
The authors are grateful for informative discussions with Patrick K\"urschner
(HTWK Leipzig), Siobhan Correnty and Elias Jarlebring (both KTH Stockholm), and Omar
De La Cruz and Lothar Reichel (both Kent State University). We thank
Omar De La Cruz for pointing out \cite{alghamdi2022greedy} and
\cite{ruymbeek2022tensor}.

\section*{Availability of data and materials}

The code used for the numerical experiments is available from GitHub,
\url{https://github.com/thomasmach/PEVP_with_Taylor_and_Chebyshev}.

\end{document}

%% file: code/fig41.tex
\begin{figure}
\centering
\includegraphics{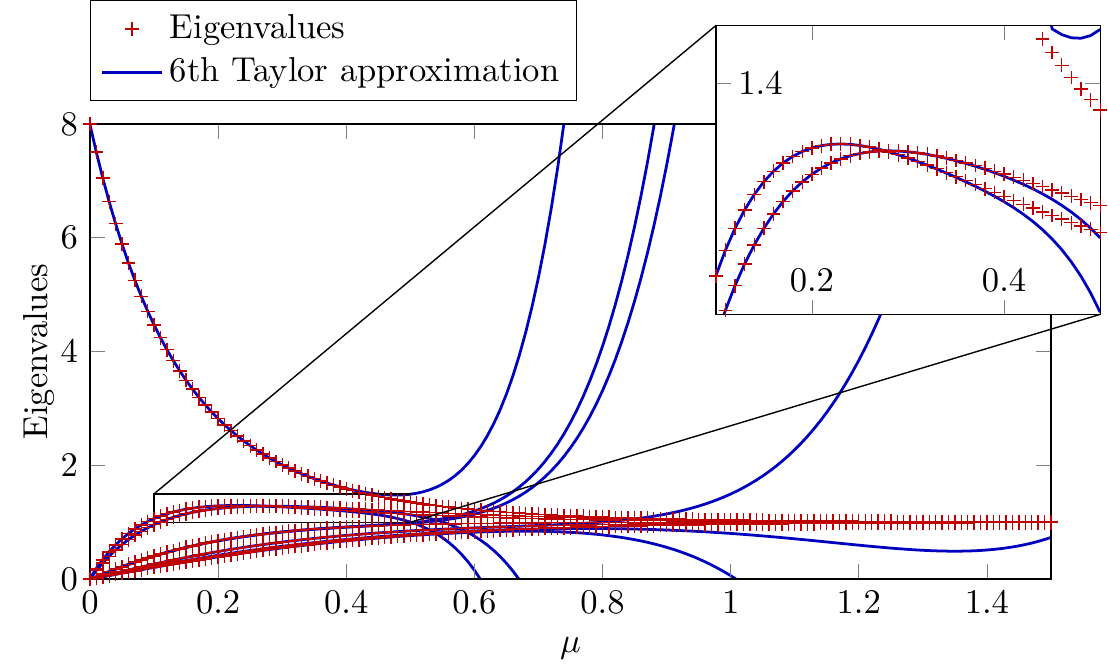}
\caption{Example~\ref{example:1}, $n=8$, $\mu_{0}=0.20$.}%
\label{fig:example1_1_8_7_0.20_1}%
\end{figure}

%% file: code/fig52.tex
\begin{figure}
\centering
\includegraphics{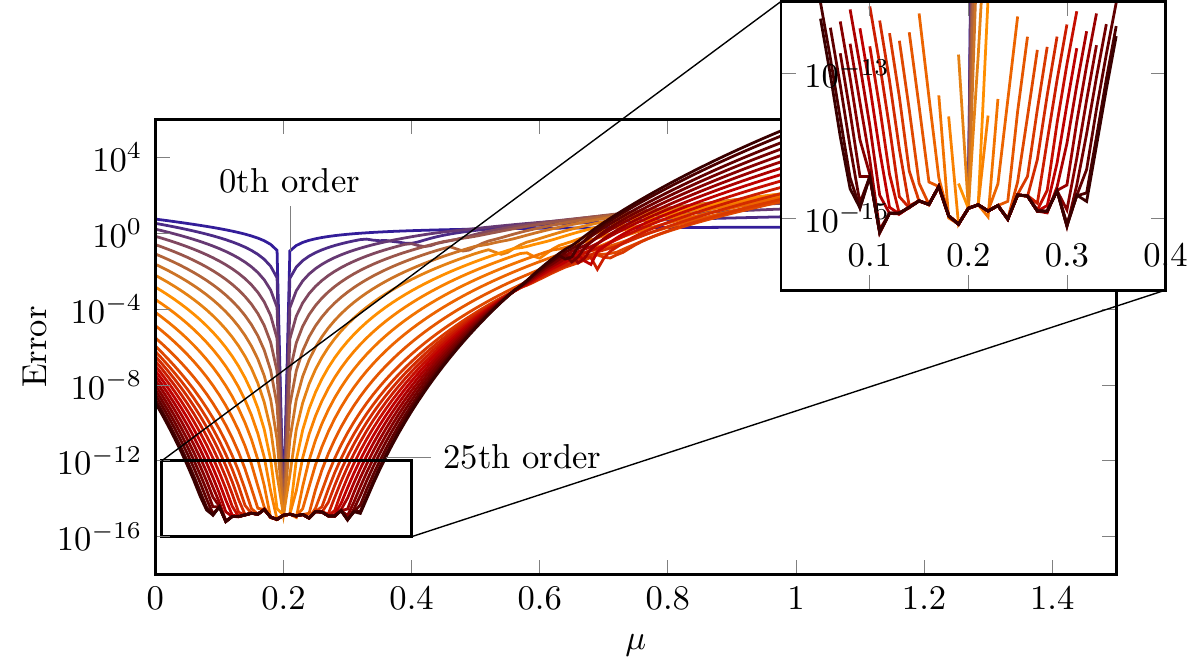}
\caption{Example~\ref{example:1}---Maximal difference between Taylor approximations and eigenvalues, $\mu_{0}=0.20$.}%
\label{fig:example4_1_8_26_0.20_1}%
\end{figure}

%% file: code/fig53.tex
\begin{figure}
\centering
\includegraphics{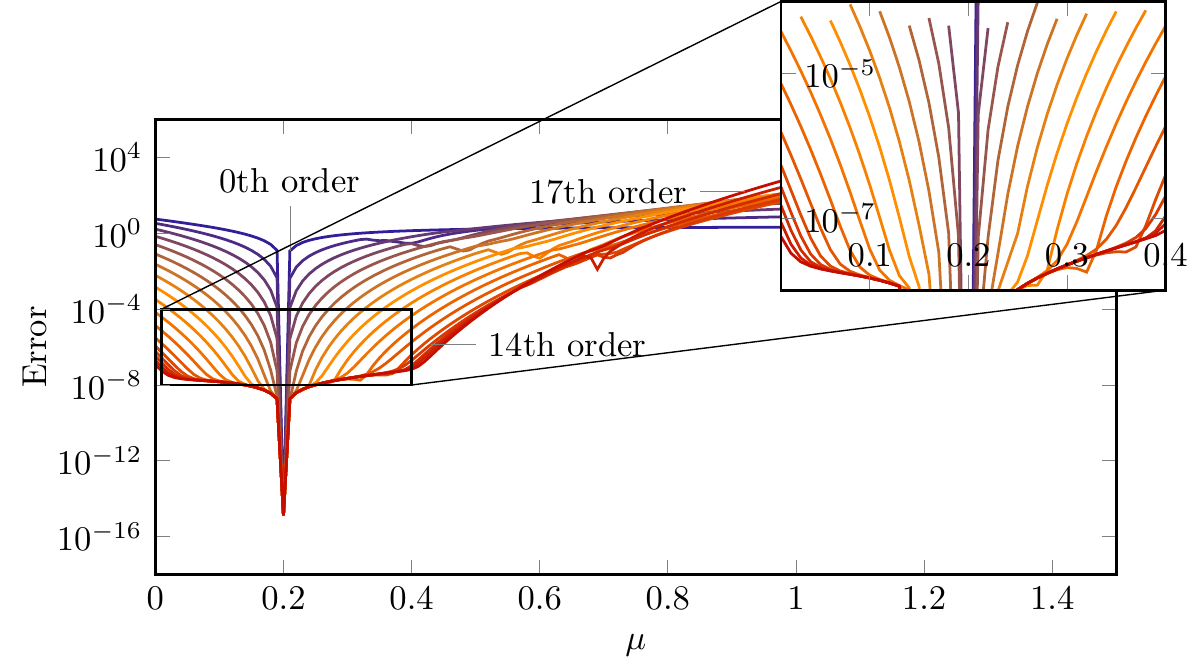}
\caption{Example~\ref{example:1}---Maximal difference between Taylor approximations and eigenvalues using single precision to store $E$, $\mu_{0}=0.20$.}%
\label{fig:example4vpa_1_8_18_0.20_1}%
\end{figure}

%% file: code/fig52_cheb.tex
\begin{figure}
\centering
\includegraphics{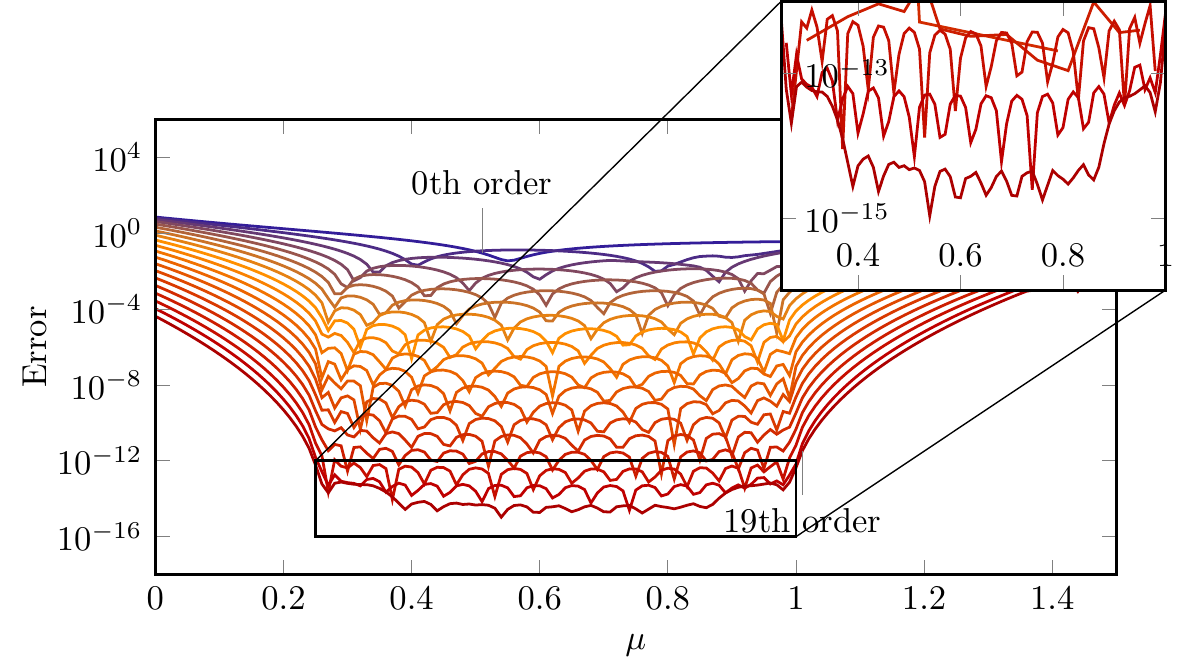}
\caption{Example~\ref{example:1}---Maximal difference between Chebyshev approximations and eigenvalues, $[\mu_{1},\mu_{2}] = [0.25,1.00]$.}%
\label{fig:example4_0_8_20_0.25_1}%
\end{figure}

%% file: code/fig61_both.tex
\begin{figure}
\centering
\includegraphics{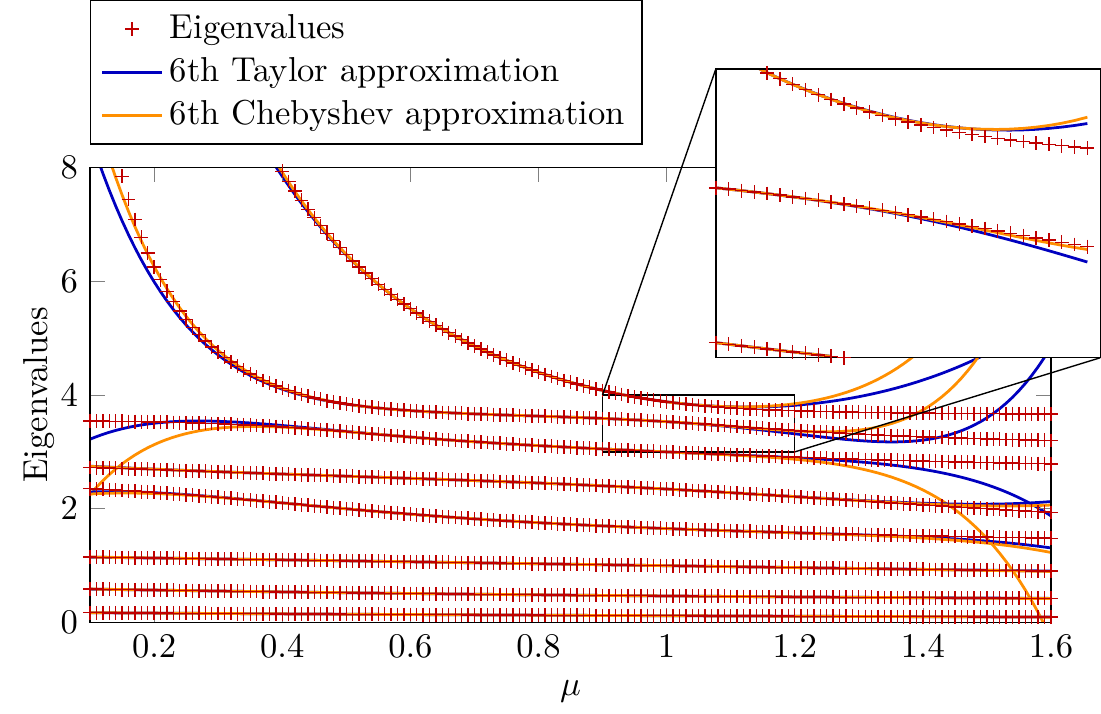}
\caption{Example~\ref{example:2}, $n=8$, $\mu_{0}=0.80$ and $[\mu_{1},\mu_{2}]=[0.50,1.00]$, respectively.}%
\label{fig:example1_1_8_7_0.80_2}%
\end{figure}

%% file: code/fig62.tex
\begin{figure}
\centering
\includegraphics{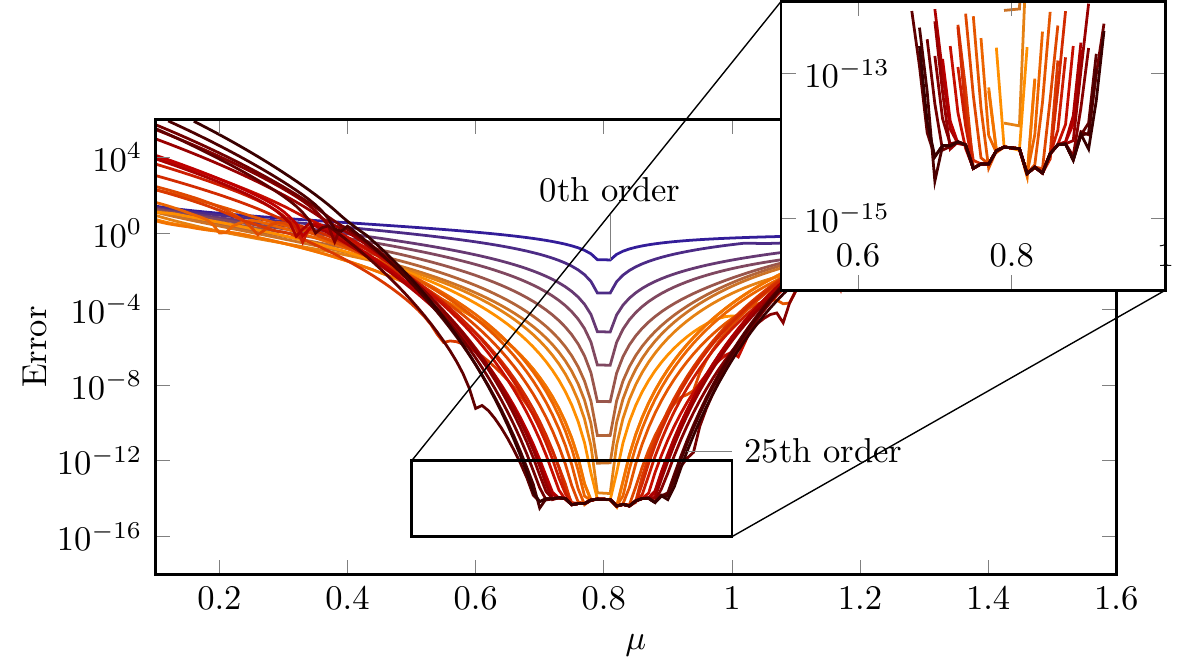}
\caption{Example~\ref{example:2}---Maximal difference between Taylor approximations and eigenvalues, $\mu_{0}=0.80$.}%
\label{fig:example4_1_8_26_0.80_2}%
\end{figure}

%% file: code/fig62_cheb.tex
\begin{figure}
\centering
\includegraphics{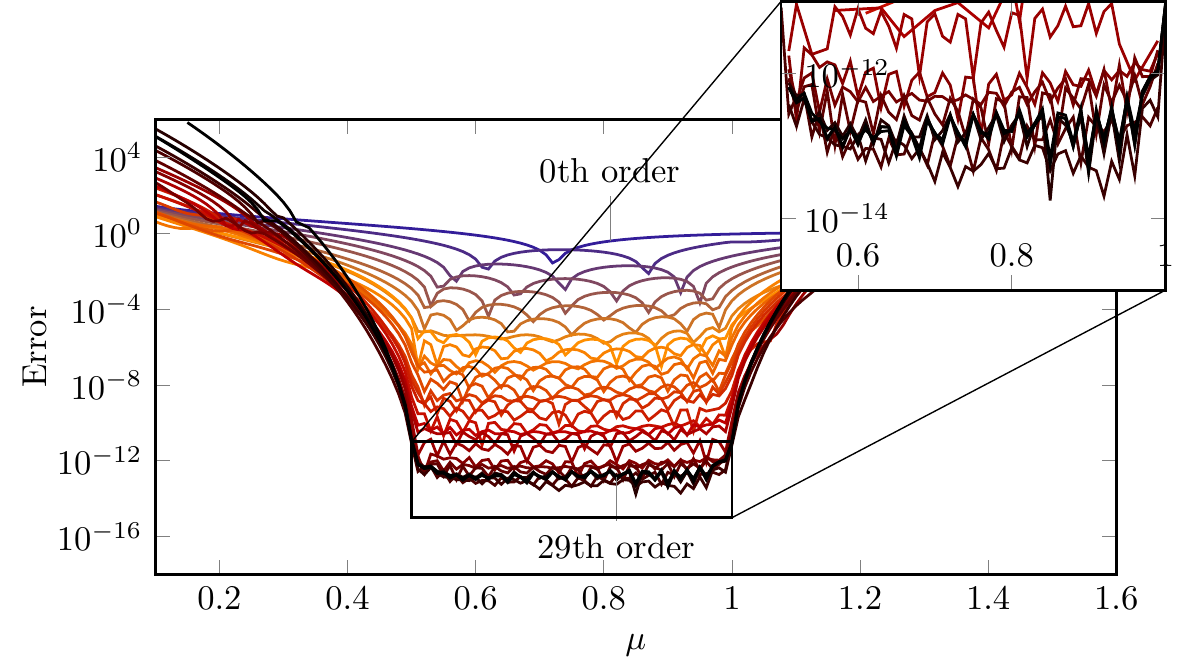}
\caption{Example~\ref{example:2}---Maximal difference between Chebyshev approximations and eigenvalues, $[\mu_{1},\mu_{2}] = [0.50,1.00]$.}%
\label{fig:example4_0_8_30_0.50_2}%
\end{figure}

%% file: code/fig71.tex
\begin{figure}
\centering
\includegraphics{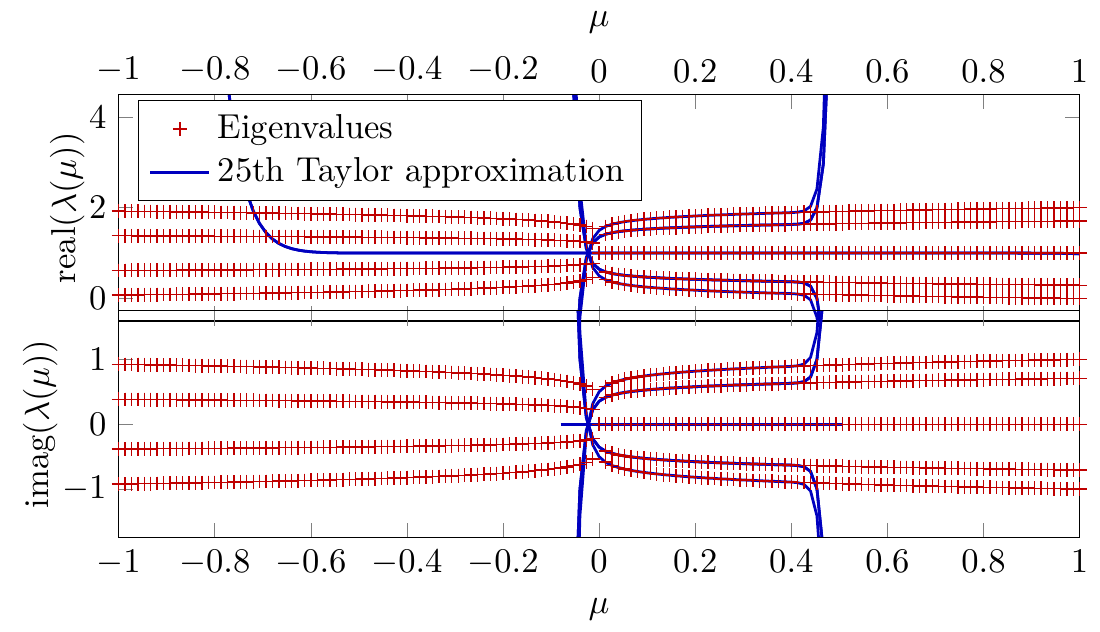}
\caption{Example~\ref{example:3}, $n=8$, $\mu_{0}=0.20$.}%
\label{fig:example1_1_8_26_0.20_3}%
\end{figure}

%% file: code/fig71_cheb.tex
\begin{figure}
\centering
\includegraphics{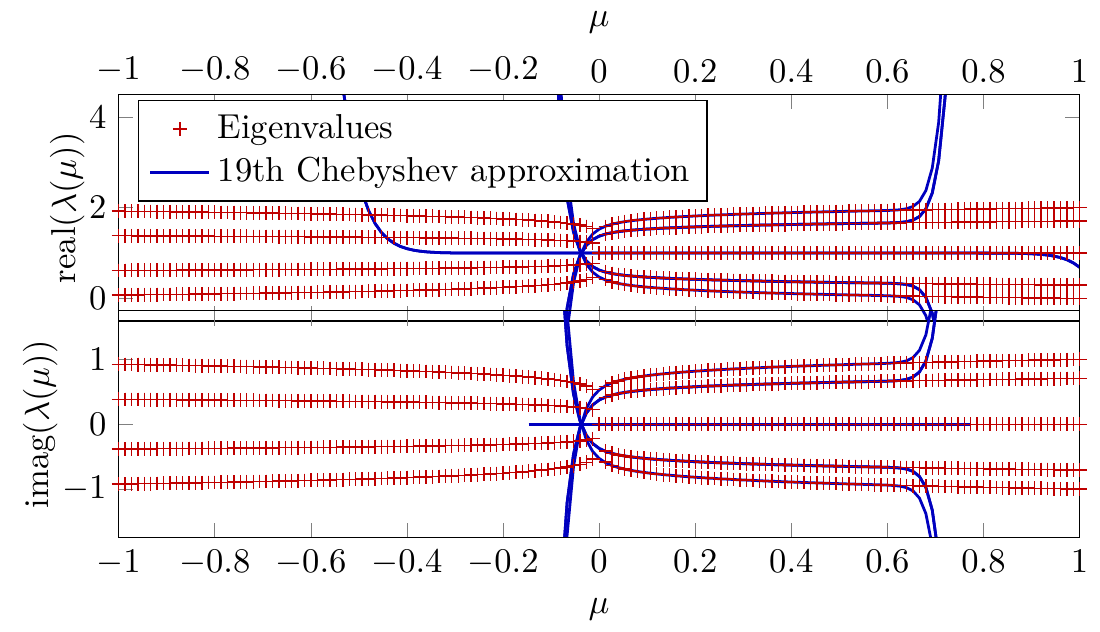}
\caption{Example~\ref{example:3}, $n=8$, $[\mu_{1},\mu_{2}]=[0.10,0.50]$.}%
\label{fig:example1_0_8_20_0.10_3}%
\end{figure}

%% file: code/fig61_taylor_1.tex
\begin{figure}
\centering
\includegraphics{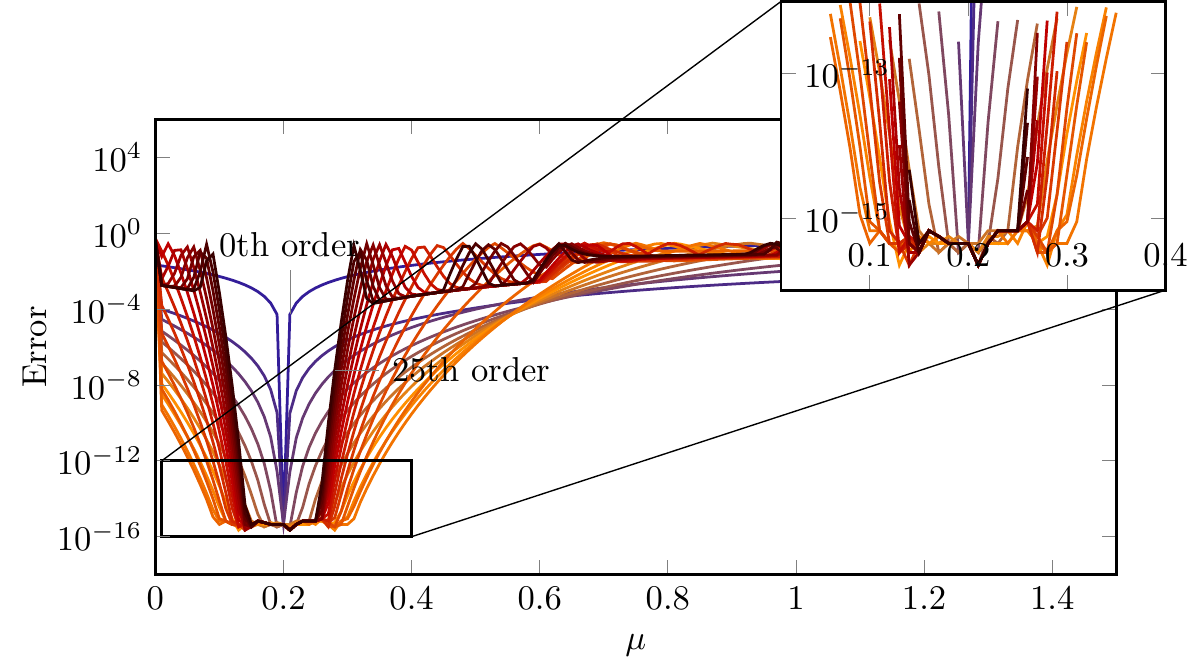}
\caption{Example~\ref{example:1}---Maximal difference between Taylor approximations and eigenvectors, $\mu_{0}=0.20$.}%
\label{fig:example5_1_8_26_0.20_1}%
\end{figure}

%% file: code/fig61_taylor_2.tex
\begin{figure}
\centering
\includegraphics{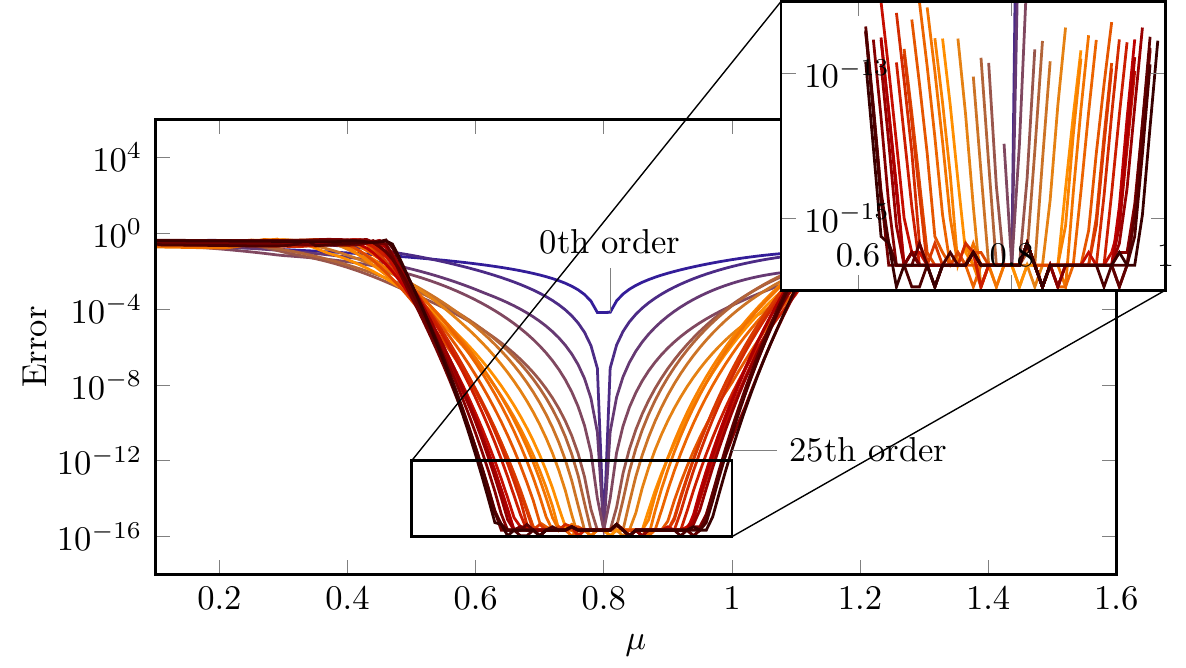}
\caption{Example~\ref{example:2}---Maximal difference between Taylor approximations and eigenvectors, $\mu_{0}=0.80$.}%
\label{fig:example5_1_8_26_0.80_2}%
\end{figure}

%% file: code/fig61_cheb_1.tex
\begin{figure}
\centering
\includegraphics{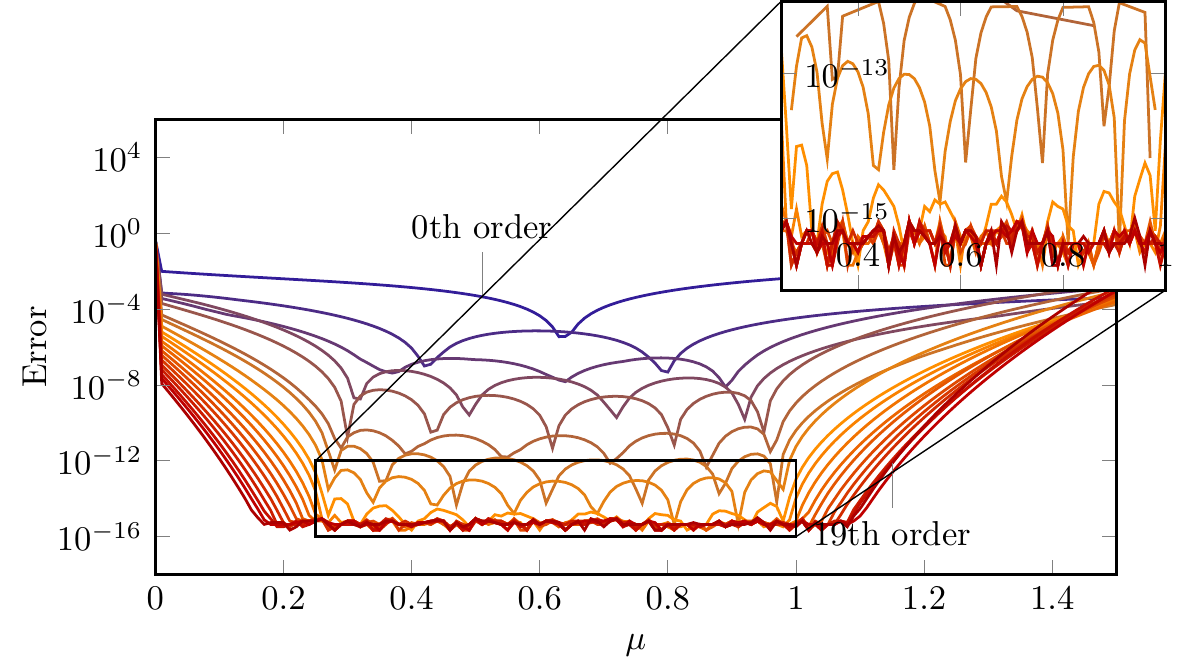}
\caption{Example~\ref{example:1}---Maximal difference between Chebyshev approximations and eigenvectors, $[\mu_{1},\mu_{2}] = [0.25,1.00]$.}%
\label{fig:example5_0_8_20_0.25_1}%
\end{figure}

%% file: code/fig61_cheb_2.tex
\begin{figure}
\centering
\includegraphics{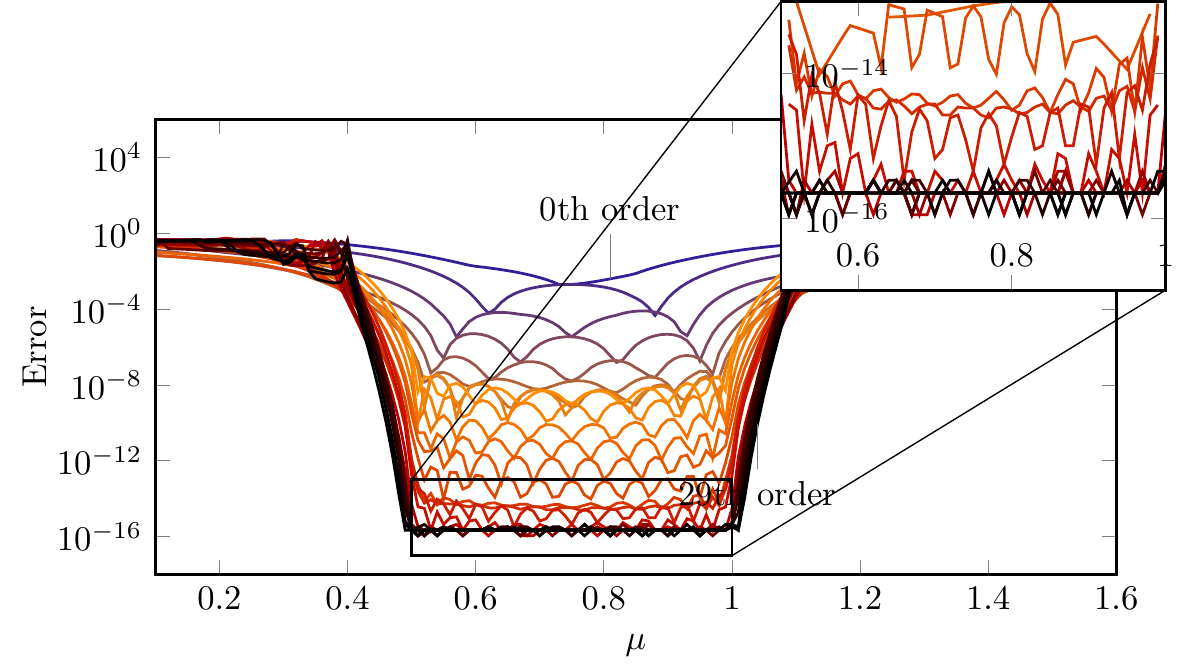}
\caption{Example~\ref{example:2}---Maximal difference between Chebyshev approximations and eigenvectors, $[\mu_{1},\mu_{2}] = [0.50,1.00]$.}%
\label{fig:example5_0_8_30_0.50_2}%
\end{figure}

%% file: code/fig71_cheb_1.tex
\begin{figure}
\centering
\includegraphics{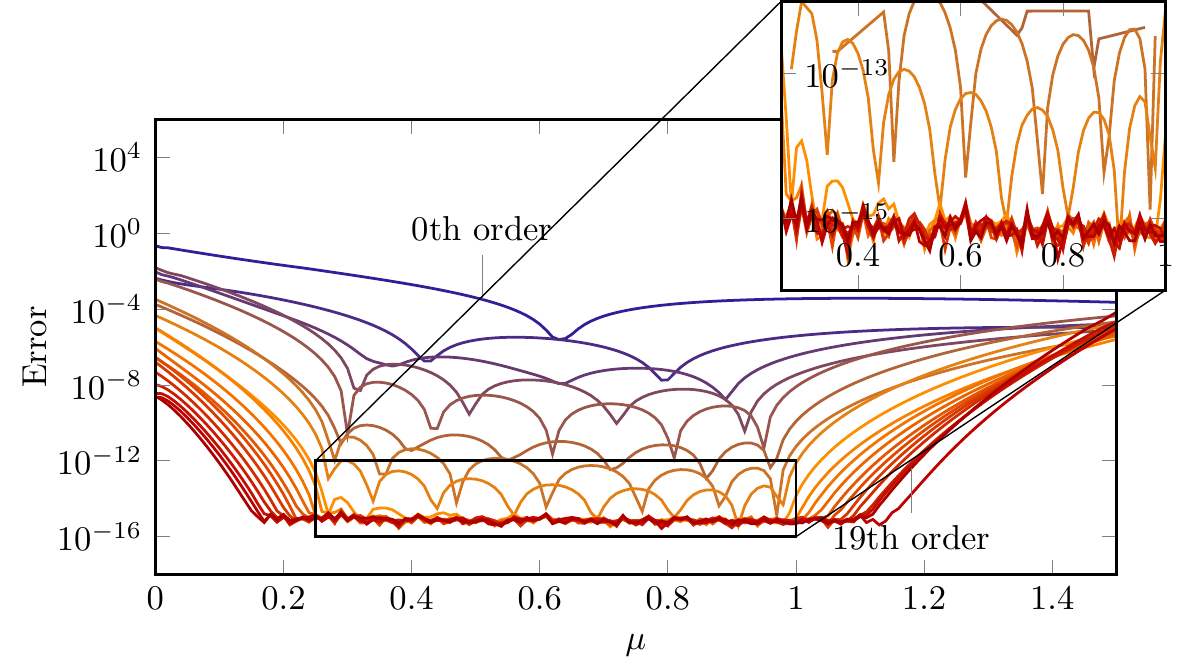}
\caption{Example~\ref{example:1}---Maximal difference between Rayleigh quotient of the Chebyshev eigenvector approximations and the eigenvalues, $[\mu_{1},\mu_{2}] = [0.25,1.00]$.}%
\label{fig:example7_0_8_20_0.25_1}%
\end{figure}